\documentclass[preprint,12pt,authoryear]{elsarticle}

\usepackage{amssymb}
\usepackage[margin=0.9in]{geometry}

\usepackage{graphics}
\usepackage{amsthm}
\usepackage{amsmath,amssymb,amsfonts}
\usepackage{mathrsfs}
\usepackage{wasysym}
\usepackage{algorithm}
\usepackage{algorithmic}
\usepackage{enumitem}
\usepackage{bm}
\usepackage{diagbox}

\usepackage{natbib}
\setcitestyle{numbers,square}
\usepackage{multirow,booktabs,comment}
\usepackage{epstopdf}

\usepackage{todonotes}
\usepackage[normalem]{ulem}

\usepackage{xcolor}
\usepackage{soul}

\usepackage[colorlinks=true,urlcolor=blue,citecolor=red,
anchorcolor=black,linkcolor=blue]{hyperref}

\usepackage{subcaption}
\usepackage[toc,page]{appendix}

\usetikzlibrary{calc,trees,positioning,arrows,chains,shapes.geometric,%
	decorations.pathreplacing,decorations.pathmorphing,shapes,%
	matrix,shapes.symbols}

\tikzset{
	>=stealth',
	punktchain/.style={
		rectangle, 
		rounded corners, 
		draw=black, very thick,
		text width=10em, 
		minimum height=2em, 
		text centered, 
		on chain},
	line/.style={draw, thick, <-},
	element/.style={
		tape,
		top color=white,
		bottom color=blue!50!black!60!,
		minimum width=6em,
		draw=blue!40!black!90, very thick,
		text width=10em, 
		minimum height=3.5em, 
		text centered, 
		on chain},
	every join/.style={->, thick,shorten >=1pt},
	decoration={brace},
	tuborg/.style={decorate},
	tubnode/.style={midway, right=2pt},
}

\tikzstyle{arrow} = [thick,->,>=stealth]

\newtheorem{thm}{Theorem}

\newtheorem{assump}{Assumption}

\newcommand{\bx}{\boldsymbol{x}}

\newcommand{\btheta}{\boldsymbol{\theta}}

\newcommand{\norm}[2]{\left\| #1 \right\|_{#2}}

\newcommand{\mb}{\boldsymbol}

\newcommand{\xs}{\mathbb}

\newcommand{\comm}[1]{}

\DeclareMathOperator*{\argmin}{arg\,min}

\usepackage{lineno}

\journal{}

\makeatletter
\def\ps@pprintTitle{%
	\let\@oddhead\@empty
	\let\@evenhead\@empty
	\let\@oddfoot\@empty
	\let\@evenfoot\@oddfoot
}
\makeatother

\begin{document}

\begin{frontmatter}



\title{Deep adaptive sampling for surrogate modeling without labeled data}

\author[mymainaddress]{Xili Wang}
\ead{xiliwang@stu.pku.edu.cn}

\author[mysecondaddress]{Kejun Tang\corref{mycorrespondingauthor}}
\cortext[mycorrespondingauthor]{Corresponding author}
\ead{tangkejun@icode.pku.edu.cn}

\author[mythirdaddress]{Jiayu Zhai}
\ead{zhaijy@shanghaitech.edu.cn}

\author[myfourthaddress]{Xiaoliang Wan}
\ead{xlwan@lsu.edu}

\author[mymainaddress,mysecondaddress]{Chao Yang}
\ead{chao\_yang@pku.edu.cn}

\address[mymainaddress]{School of Mathematical Sciences, Peking University}

\address[mysecondaddress]{PKU-Changsha Institute for Computing and Digital Economy}

\address[mythirdaddress]{Institute of Mathematical Sciences, ShanghaiTech University }

\address[myfourthaddress]{Department of Mathematics and Center for Computation and Technology, Louisiana State University}

\begin{abstract}
        Surrogate modeling is of great practical significance for parametric differential equation systems. In contrast to classical numerical methods,  
        using physics-informed deep learning methods to construct simulators for such systems is a promising direction due to its potential to handle 
        high dimensionality, which requires minimizing a loss over a training set of random samples. However, the random samples introduce statistical errors, which may become the dominant errors for the approximation of low-regularity and high-dimensional problems. In this work, we present a deep adaptive sampling method for surrogate modeling ($\mathrm{DAS}^2$), where we generalize the deep adaptive sampling (DAS) method \citep{tang_das} [Tang, Wan and Yang, 2023] to build surrogate models for low-regularity parametric differential equations. In the parametric setting, the residual loss function 
        can be regarded as an unnormalized probability density function (PDF) of the spatial and parametric variables. This PDF is approximated by a deep generative model, from which new samples are generated and added to the training set. Since the new samples match the residual-induced distribution, the refined training set can further reduce the statistical error in the current approximate solution. 
        We demonstrate the effectiveness of $\mathrm{DAS}^2$ with a series of numerical experiments, including the parametric lid-driven 2D cavity flow problem with a continuous range of Reynolds numbers from 100 to 1000.
\end{abstract}

%


\end{frontmatter}


\section{Introduction}
Solving differential equations with different parametric settings is widely found in uncertainty quantification \cite{xiu2002wiener,xiu2010numerical,xiu2016stochastic}, inverse design \cite{JagtKarn20,ghosh2022inverse}, Bayesian inverse problems \cite{stuart2010inverse,li2014adaptive,cui2016scalable,liao2019adaptive,xia2022bayesian, feng2023dimension}, digital twins \cite{chakraborty2021role,kapteyn2022digitaltwins,torzoni2024digital}, parametric optimal control \cite{yin2023aonn}, and shape optimization \cite{wang2023aonn}, etc. The computational cost of solving such parametric differential equations with conventional numerical methods is expensive because repeated simulations (i.e., \emph{many-query}) of differential equations are required. To handle such many-query problems, one may construct a surrogate model that can efficiently predict the parametric solution without sacrificing much accuracy, 
which is sufficient for many engineering applications. For instance,  
reduced order models (ROM) \cite{boyaval2010reduced,elman2013reduced,quarteroni2015reduced,chen2019robust} are widely used in practice, 
where the approximate solution is expressed as a linear combination of some bases that are computed by low-rank approximation of \emph{snapshot matrices}. 
ROM becomes inefficient if the parametric solution does not lie in a low-dimensional linear subspace 
\cite{chaturantabut2010nonlinear,bonito2021nonlinear,cohen2023nonlinear}. 

Deep learning-based methods for surrogate modeling have been proposed to give an alternative approach. One straightforward data-driven approach is to utilize deep neural networks to learn a mapping from a parameterized function space to the solution space \cite{zhu2018bayesian,lu2021learning, li2021fourier}, where simulation-based input-output pair data are used to train the deep neural networks. 
Constructing surrogate models without labeled data is necessary to handle the cases where simulation or experimental data are scarce or not available. To this end, numerical strategies have been developed for the neural network approximation of deterministic partial differential equations (PDEs) \cite{sirignano2018dgm,raissi2019physics,weinan2018deep,han2018solving, karniadakis2021physics,mao2020physics, hao2022physics}, based on which parametric PDEs can also be addressed.  
For example, physics-informed deep learning is used to construct surrogate models for efficient uncertainty quantification \cite{zhu2019physics,sun2020surrogate,WangGao21}. Depending on the formulation of parametric (partial or ordinary) differential equations, two types of neural network models can be considered: one is a plain neural network whose inputs include both the parameters and the spatial variables, and the other one is defined by an operator learning problem, e.g., DeepONet \cite{ lu2021learning,wang2021learning}. Regardless of the structure of the surrogate model, the underlying training procedure is similar, which minimizes a loss functional discretized by the random collocation points in the physical domain and the parametric space. For low-regularity deterministic problems, the collocation points significantly affect the generalization error of neural networks \cite{tang2022adaptive,tang_das}. This issue becomes worse for surrogate modeling because the parametric space may also introduce low regularity other than the additional dimensions. For example, in the Navier-Stokes equations, large Reynolds numbers cause some small-scale structures. To capture these features, the distribution of the collocation points in the training set must be consistent with the characteristics of the velocity field, and such a correspondence needs to be maintained for all Reynolds numbers considered as the inputs of a surrogate model. Hence, we need to pay particular attention to the random collocation points in the training set to obtain a sufficiently accurate surrogate model.


In this work, we develop a deep adaptive sampling approach for surrogate modeling ($\mathrm{DAS}^2$) without labeled data, which generalizes DAS \cite{tang_das} to the parametric setting. Without losing generality, neural networks with augmented parameter inputs are used to approximate the parametric solutions. We intend to find a certain set of collocation points that results in a relatively flat residual profile. Since a flat residual profile has a small variance, the statistical error in the discretization of the loss functional can be significantly reduced for a fixed number of samples, which eventually improves the accuracy of the approximate solution. The desired training set is achieved through iterations. Assume that the surrogate model is trained with respect to a certain training set. The total residual of the parametric equations is viewed as an unnormalized probability density function (PDF). More samples will be introduced to the training set in the region of high density such that the residual over there can be reduced. To achieve this, a deep generative model is trained to approximate the residual-induced PDF, and new samples are drawn from this trained deep generative model. Once the training set is updated, the surrogate model will be further trained, after which the aforementioned procedure is repeated. The same algorithm can be applied to other types of surrogate models such as DeepONet. The main contributions of this work are summarized as follows.


\subsection{Main contributions}
\begin{itemize}
	\item We propose a deep adaptive sampling approach for surrogate modeling of parametric differential equations without labeled data.
	\item We demonstrate the efficiency of the proposed method with a series of numerical experiments, including the operator learning problem, the parametric optimal control problem, and the lid-driven 2D cavity flow problem with a continuous range of Reynolds numbers from 100 to 1000.
\end{itemize}

\subsection{Related work}
The adaptive sampling-based neural network methods for solving deterministic differential equations are under active development. Nevertheless, adaptive sampling of parametric differential equations is still to be studied. We summarize the most related lines of this work: adaptive sampling methods for deterministic problems and neural network methods for parametric differential equations.  

\subsubsection{Adaptive sampling methods}\label{sec_intro_adap}
Solving (partial) differential equations with deep learning methods usually needs a large set of collocation points, particularly when the solution has subtle structures such as high frequency, high-density concentration, multiscale structure, or discontinuity \cite{olga2020limitations, raissi2020hidden, wang2022when, zhai2022adeep}. Adaptive collocation points may significantly reduce the computational cost, where the essence is to define a proper error indicator and generate training collocation points accordingly.


The residual-based adaptive refinement (RAR) method \cite{lu2021deepxde, wu2023comprehensive} is proposed to enhance the performance of physics-informed machine learning.
In RAR, one needs to construct a set of uniform samples as a candidate set, within which the samples associated with large residuals are selected and added to the current training set. 
However, 
such a strategy is not effective for high-dimensional problems since most of the volume of the computational domain concentrates around its surface \cite{wright2021high}. To obtain true samples from the residual-induced distribution, classical sampling methods such as MCMC can be employed \cite{gao2021active,yu2023mcmcpinn,wen2023coupling}, which, however, are also affected by the curse of dimensionality. To handle high-dimensional problems, we need to introduce other techniques. In \cite{tang_das}, the deep adaptive sampling (DAS) method is proposed, where a normalizing flow model is used to approximate the residual-induced distribution, based on which new collocation points are generated to further improve the accuracy of the current approximate solution. 
DAS uses the current residual as an explicit guidance for the selection of new collocation points, which is similar to the procedure of classical adaptive methods such as the adaptive finite element method. 
Another track is to implicitly search for a distribution that generates collocation points that result in a smooth residual profile.  
In \cite{tang2023adversarial}, an adversarial adaptive sampling framework (AAS) is proposed to seek an optimal model for the solution and an optimal distribution for the training set at the same time through a min-max formulation, which can be regarded as a generalization of the strategies that aim to find a better weight for each fixed sample \cite{han2022residual,Sokratis2023_weight}.
To reduce the training cost from the deep generative model, one can replace the deep generative model in DAS or AAS with other density models, such as Gaussian mixture models \cite{jiao2023gas}. In \cite{daw2022mitigating, gao2023failure2, gao2023failure}, the authors reformulate the adaptive sampling procedure as a failure event subject to a threshold that helps determine where new collocation points are needed. 
Other related works include \cite{peng2022rang,subramanian2022adaptive,chen2023adaptive,hou2023enhancing}. 


\subsubsection{Neural network methods for parametric differential equations} 
The study of parametric PDEs with neural network methods started from the very beginning of PINNs in \cite{raissi2019physics}, where neural network provides a general model for both forward and inverse problems. We focus on parametric forward problems in this work, which can also be regarded as an operator learning problem. 
In \cite{lu2021learning}, DeepONet is proposed, 
which formulates an operator that maps infinite-dimensional data, e.g., boundary and initial conditions, to the solution functions of parametric differential equations. Another 
typical operator model is the Fourier Neural Operaotr (FNO) \cite{li2020Neural, li2021fourier}. 
Physics-informed operator learning has been developed in \cite{wang2021learning,li2021physics} to reduce the dependence of DeepONet and FNO on data in the training process. 



\section{Problem setting and statistical errors in physics-informed surrogate modeling}\label{sec_probset}
Let $\Omega_s$ be a spatial domain (in $\xs{R}^n$) which is bounded, connected and with a smooth boundary $\partial \Omega_s$,
and $\mb{x} \in \Omega_s$ is a spatial variable.
Let $\mb{\xi}$ be a vector which collects a finite number of parameters. The dimension of $\mb{\xi}$ is denoted by $d$, i.e., we write $\mb{\xi}=[\xi_1,\ldots,\xi_d]^\mathsf{T}$. 
We restrict our attention to the situation that $\mb{\xi}$ has a bounded and connected support.
Without loss of generality, we next assume the support of $\mb{\xi}$ to be $\Omega_p$ where $\Omega_p: = [C_1, C_2]^d$ and $C_1, C_2$ are two constants.
The physics of problems considered here
are governed by 
differential equations over the spatial domain $\Omega_s$ and
boundary conditions on the boundary $\partial \Omega_s$.
Consider the following parametric differential equations: find $u: \Omega_s \times \Omega_p \mapsto \xs{R}$ such that
\begin{align}
	\mathcal{L}\left(\mb{x},\mb{\xi};u\left(\mb{x}, \mb{\xi} \right)\right)=s(\mb{x}, \mb{\xi})  \qquad
	&\forall \left(\mb{x},\mb{\xi} \right) \in \Omega_s \times \Omega_p,
	\label{spdexi1}\\
	\mathcal{B}\left(\mb{x},\mb{\xi};u\left(\mb{x},\mb{\xi} \right)\right)=g(\mb{x}, \mb{\xi}) \qquad
	&\forall \left(\mb{x}, \mb{\xi} \right)\in \partial \Omega_s \times \Omega_p,
	\label{spdexi2}
\end{align}
where $\mathcal{L}$ is a differential operator and $\mathcal{B}$ is a boundary
operator, both of which can involve parameters. $s$ is the source function and $g$ specifies the boundary conditions. To simplify the notation, we denote $\Omega = \Omega_s \times \Omega_p$ and $\partial \Omega = \partial \Omega_s \times \Omega_p$. The goal of this study is to construct a surrogate model, which is the parameterized solution to the differential equation, by using the deep adaptive sampling method (DAS) \cite{tang_das} without labeled data, i.e., DAS for surrogate modeling, and we call this strategy $\mathrm{DAS}^2$ for short. Once this surrogate is constructed, the solution $u(\mb{x},\mb{\xi})$ can be efficiently predicted for any $\mb{\xi}$ without solving the (partial or ordinary) differential equation repeatedly. 

The framework of physics-informed surrogate modeling for parametric differential equations is as follows. Let $u_{\mb{\theta}}(\mb{x}, \mb{\xi})$ be a neural network parameterized with $\mb{\theta}$, where the input of the neural network is the tuple $(\mb{x}, \mb{\xi})$. One can use $u_{\mb{\theta}}(\mb{x}, \mb{\xi})$ to approximate $u(\mb{x},\mb{\xi})$ through minimizing the following loss functional
\begin{equation}\label{eq_s_loss}
	\begin{aligned}
		J \left( u_{\mb{\theta}} \right) 
		&= J_r(u_{\mb{\theta}}) + \gamma J_b(u_{\mb{\theta}}) \quad \text{with} \\
		J_r(u_{\mb{\theta}}) &= \int_{\Omega} |r(\mb{x},\mb{\xi};\mb{\theta})|^2 d\mb{x}d\mb{\xi} \ \text{ and } \ 
		J_b(u_{\mb{\theta}})= \int_{\partial \Omega} |b(\mb{x},\mb{\xi};\mb{\theta})|^2 d\mb{x}d\mb{\xi},
	\end{aligned}
\end{equation} 
where $r(\mb{x},\mb{\xi};\mb{\theta}) = \mathcal{L} u_{\mb{\theta}}(\mb{x},\mb{\xi}) - s(\mb{x},\mb{\xi})$, and $b(\mb{x},\mb{\xi};\mb{\theta}) = \mathcal{B} u_{\mb{\theta}}(\mb{x},\mb{\xi}) - g(\mb{x},\mb{\xi})$ are the residuals that measure how well $u_{\mb{\theta}}$ satisfies the parametric differential equations and the boundary conditions, respectively, and $\gamma>0$ is a penalty parameter. Before optimizing this loss functional with respect to $\mb{\theta}$, we need to discretize the integral defined in \eqref{eq_s_loss} numerically. In general, choosing uniformly distributed collocation points is a standard way for discretizing the integral.
Let $\mathsf{S}_{\Omega} = \{\mb{x}_{\Omega}^{(i)},\mb{\xi}^{(i)} \}_{i=1}^{N_r}$ and $\mathsf{S}_{\partial \Omega} = \{\mb{x}_{\partial \Omega}^{(i)}, \mb{\xi}^{(i)} \}_{i=1}^{N_b}$ be two sets of uniformly distributed collocation points respectively. 
We then minimize the following empirical loss in practice
\begin{equation}\label{eq_discrete_residual}
	J_N \left( u_{\mb{\theta}} \right) = J_{r,N} + \gamma J_{b,N} = \frac{1}{N_r} \sum\limits_{i=1}^{N_r} r^2(\mb{x}_{\Omega}^{(i)},\mb{\xi}^{(i)};\mb{\theta}) + \gamma \frac{1}{N_b} \sum\limits_{i=1}^{N_b} b^2(\mb{x}_{\partial \Omega}^{(i)},\mb{\xi}^{(i)};\mb{\theta}),
\end{equation}
which can be regarded as the Monte Carlo (MC) approximation of $J(u_{\mb{\theta}})$ subject to a statistical error of $\mathit{O}(N^{-1/2})$ with $N$ being the sample size. Let $u_{\mb{\theta}_N^*}$ be the minimizer of the empirical loss $J_N(u_{\mb{\theta}})$
\begin{equation}\label{empminimizer}
	u_{\mb{\theta}_N^*} = \arg \min_{\mb{\theta}} J_N(u_{\mb{\theta}})
\end{equation}
and $u_{\mb{\theta}^*}$ be the minimizer of the original loss functional $J(u_{\mb{\theta}})$
\begin{equation}\label{trueminimizer}
	u_{\mb{\theta}^*} = \arg \min_{\mb{\theta}} J(u_{\mb{\theta}}).
\end{equation}
We can decompose the error of $u_{\mb{\theta}_N^*}$ into two parts as follows
\begin{equation*}
	\mathbb{E} \left( \norm{u_{\mb{\theta}_N^*} - u}{\Omega} \right) \leq \mathbb{E} \left( \norm{u_{\mb{\theta}_N^*} - u_{\mb{\theta}^*}}{\Omega} \right) + \norm{u_{\mb{\theta}^*} - u}{\Omega},
\end{equation*}
where $\mathbb{E}$ denotes the expectation with respect to the random samples and the norm $\norm{\cdot}{\Omega}$ corresponds to a function space for $u$. 
Without taking into account the optimization error, 
one can see that the total error of neural network approximation for parametric differential equations mainly consists of two parts: the approximation error and the statistical error. The approximation error is dependent on the model capability of neural networks, while the statistical error originates from the collocation points.

Uniformly distributed collocation points are not effective for training neural-network-based surrogate models if the solution has low regularity \cite{tang2022adaptive, tang_das, wu2023comprehensive} since the effective sample size of the Monte Carlo approximation of $J(u_{\mb{\theta}})$ is significantly reduced by the large variance induced by the low regularity. For high-dimensional problems, random samples becomes more localized due to the curse of dimensionality \cite{wright2021high}, which shares some similarities with the low-dimensional problems of low regularity. Therefore, adaptive sampling is needed. In this work, we propose a deep adaptive sampling approach for surrogate modeling of parametric differential equations without labeled data, which generalizes the DAS method \cite{tang_das} to parametric settings. For simplicity and clarity, we only consider $J_r(u_{\mb{\theta}})$ and remove the boundary term $J_b(u_{\mb{\theta}})$. This is because one can employ some penalty-free techniques \cite{berg2018unified, sheng2020pfnn} to remove $J_b(u_{\mb{\theta}})$ from the loss. 

\section{Deep adaptive sampling for surrogate modeling}
The statistical error comes from the discretization of $J_N(u_{\mb{\theta}})$. One straightforward way to reduce the error of $J_{N}(u_{\mb{\theta}})$ is to increase the number of uniformly distributed collocation points in the training set. However, if the solution is of low regularity, the large variance of the residual will significantly reduce the number of effective samples for the computation of $J_N(u_{\mb{\theta}})\approx J(u_{\mb{\theta}})$ such that the final approximate solution may gain barely any improvement. To alleviate this issue, the selection of collocation points must be consistent with the problem properties, in other words, adaptive sampling needs to be considered. 

A deep adaptive sampling (DAS) method has been developed in \cite{tang_das} for deterministic PDEs. We in this work intend to generalize DAS to deal with parametric differential equations
and call this generalization $\mathrm{DAS}^2$ for short, i.e., \emph{deep adaptive sampling} for \emph{surrogates}. The main difficulties come from the additional dimensions from $\mb{\xi}$. First, samples are needed from both the spatial domain and the parametric space for the discretization of the loss functional. 
Second, low regularity may come from the spatial domain, the parametric space, or both. Without assuming any prior knowledge of the residual profile, we need to efficiently generate random samples that are consistent with an arbitrary high-dimensional distribution. To handle such a situation, we employ a deep generative model, called KRnet, to approximate the residual-induced distribution and then generate random collocation points accordingly.
The PDF defined by KRnet is 
\begin{equation}\label{eq_krpdf}
	p_{\mathsf{KRnet}}(\mb{x},\mb{\xi};\mb{\theta}_f)=p_{\mb{Z}}(f_{\mathsf{KRnet}}(\mb{x},\mb{\xi};\mb{\theta}_f)) \left |\det\nabla_{\mb{x},\mb{\xi}} f_{\mathsf{KRnet}} \right|,
\end{equation}
where $f_{\mathsf{KRnet}}$ denotes an \emph{invertible} mapping defined by KRnet parameterized with $\mb{\theta}_f$, and the prior distribution $p_{\mb{Z}}$ for the random vector $\mb{z}$ is usually chosen as the standard normal distribution. 
The overall structure of KRnet is specified as follows
\begin{equation*}
	\mb{z} = f_{\mathsf{KRnet}}(\mb{x},\mb{\xi}) = L_{N} \circ f_{[K-1]}^{\textsf{outer}} \circ \cdots \circ f_{[1]}^{\textsf{outer}} (\mb{x},\mb{\xi}),
\end{equation*}
where $f_{[i]}^{\textsf{outer}}$ is defined as
\begin{equation*}
f_{[k]}^{\textsf{outer}} = L_S \circ f_{[k, L]}^{\textsf{inner}} \circ \cdots \circ f_{[k,1]}^{\textsf{inner}} \circ L_R.
\end{equation*}
Here, $f_{[k,i]}^{\textsf{inner}}$ is a combination of $L$ affine coupling layers \cite{dinh2016density,kingma2018glow} and one scale and bias layer, and $L_N$, $L_S$ and $L_R$ represent the nonlinear layer, the squeezing layer and the rotation layer respectively, where details can be found in the literature \cite{tangwandensity2020,tang2022adaptive,tang_das,wan2021augmented}. 


\subsection{Sample from a joint PDF}
When the low regularity is related to both $\mb{x}$ and $\mb{\xi}$, the adaptive sampling for both $\mb{x}$ and $\mb{\xi}$ is needed. 
We need to generate samples from a joint PDF $\hat{r}(\mb{x},\mb{\xi})$ induced by the residual $r(\mb{x},\mb{\xi};\mb{\theta})$ for a certain $\mb{\theta}$. Following \cite{tang_das}, $\hat{r}(\mb{x},\mb{\xi})$ is defined as 
\begin{equation*}
	\hat{r}(\mb{x},\mb{\xi}) \propto r^2(\mb{x},\mb{\xi};\mb{\theta}) h(\mb{x},\mb{\xi}),
\end{equation*}
where $h(\mb{x},\mb{\xi})$ is a cutoff function as defined in \cite{tang_das}. The cutoff function $h(\mb{x},\mb{\xi})$ is defined on a compact support $B\supset\Omega$, where $h(\mb{x},\mb{\xi})=1$ if $(\mb{x},\mb{\xi})\in\Omega$ and then decays linearly to 0 towards $\partial B$. $B$ is chosen to be slightly larger than $\Omega$ \cite{tang_das}.  
We then employ the PDF model induced by KRnet to approximate $\hat{r}(\mb{x},\mb{\xi})$ on $B$.  
Mathematically, we need to solve the following optimization problem 
\begin{equation}\label{eq_kl_krres}
	\mb{\theta}_f^* = \arg \min_{\mb{\theta}_f} D_{\mathsf{KL}}(\hat{r}(\mb{x},\mb{\xi})||p_{\mathsf{KRnet}}(\mb{x},\mb{\xi};\mb{\theta}_f)),
\end{equation}
where $D_{\mathsf{KL}}(\cdot || \cdot)$ denotes the Kullback-Leibler (KL) divergence between two distributions. Let $\mb{\theta}_f^*$ be the optimal parameter. Since $B$ is slightly larger than $\Omega$, we may generate random samples as
\begin{equation*}
	(\mb{x},\mb{\xi}) = f_{\mathsf{KRnet}}^{-1}(\mb{z};\mb{\theta}_f^*),
\end{equation*}
and only keep those that belong to $\Omega$. The KL divergence in \eqref{eq_kl_krres} is 
\begin{equation*}
	D_{\mathsf{KL}}(\hat{r}(\mb{x},\mb{\xi})||p_{\mathsf{KRnet}}(\mb{x},\mb{\xi};\mb{\theta}_f)) = \int_{B} \hat{r}(\mb{x},\mb{\xi}) \log \hat{r}(\mb{x},\mb{\xi})d\mb{x}d{\mb{\xi}} - \int_{B} \hat{r}(\mb{x},\mb{\xi}) \log p_{\mathsf{KRnet}}(\mb{x},\mb{\xi})d\mb{x}d{\mb{\xi}}.
\end{equation*}
The first term is independent on $\mb{\theta}_f$, which does not affect the optimization step for $p_{\mathsf{KRnet}}$ defined in equation \eqref{eq_krpdf}. So, the PDF approximation step is equivalent to minimizing the cross entropy between $\hat{r}$ and $p_{\mathsf{KRnet}}$ \cite{de2005tutorial, rubinstein2013cross}:
\begin{equation*}
	H(\hat{r}, p_{\mathsf{KRnet}}) = -\int_{B} \hat{r} \log p_{\mathsf{KRnet}} d\mb{x}d\mb{\xi}.
\end{equation*}
To compute this cross entropy numerically, we need to use the importance sampling technique since the samples from $\hat{r}$ are not available. Here, we use a PDF model with known parameters $\hat{\mb{\theta}}_f$ for importance sampling:
\begin{equation}\label{eq_discrete_ce}
	H(\hat{r}, p_{\mathsf{KRnet}}) \approx -\frac{1}{m} \sum\limits_{i=1}^m \frac{\hat{r}(\mb{x}^{(i)},\mb{\xi}^{(i)}) \log p_{\mathsf{KRnet}}(\mb{x}^{(i)}, \mb{\xi}^{(i)};\mb{\theta}_f) }{p_{\mathsf{KRnet}}(\mb{x}^{(i)},\mb{\xi}^{(i)};\hat{\mb{\theta}}_f)},
\end{equation}
where $m$ is the number of collocation points for estimating the cross entropy and the choice of $\hat{\mb{\theta}}_f$ is specified in Algorithm \ref{alg_das_jo}. 


\subsection{Sample from a marginal PDF}
If the low regularity originates only from the parametric space, we can use a marginal PDF for adaptive sampling to reduce the complexity. 
We let
\begin{equation*}
	\tilde{r}^2(\mb{\xi};\mb{\theta}) = \int_{\Omega_s} r^2(\mb{x},\mb{\xi};\mb{\theta}) d \mb{x}. 
\end{equation*}
We assume that for any $\mb{\xi}$ and $\mb{\theta}$, $\tilde{r}^2(\mb{\xi};\mb{\theta})$ can be well approximated by a fixed set of uniform samples $\{\mb{x}^{(i)}\}_{i=1}^{m_x}$ in the spatial domain, i.e.,
\begin{equation}\label{eq_mgpdf_discrete}
	\tilde{r}^2(\mb{\xi};\mb{\theta}) \approx \frac{1}{m_{\mb{x}}}\sum\limits_{i=1}^{m_{\mb{x}}} r^2(\mb{x}^{(i)},\mb{\xi};\mb{\theta}). 
\end{equation}
In this way, the empirical loss in equation \eqref{eq_discrete_residual} could be rewritten as 
\begin{equation}\label{eq_discrete_residual_rd}
	J_N \left( u_{\mb{\theta}} \right) = \frac{1}{N_{\tilde{r}}} \sum\limits_{i=1}^{N_{\tilde{r}}} \tilde{r}^2(\mb{\xi}^{(i)};\mb{\theta}) + \gamma \frac{1}{N_b} \sum\limits_{i=1}^{N_b} b^2(\mb{x}_{\partial \Omega}^{(i)},\mb{\xi}^{(i)};\mb{\theta}).
\end{equation}
Similar to sampling a joint PDF, we can approximate the residual-induced distribution $\hat{r}(\mb{\xi})\propto \tilde{r}^2(\mb{\xi};\mb{\theta})h(\mb{\xi})$ by the following optimization problem
\begin{equation*}
	\mb{\theta}_f^* = \arg \min_{\mb{\theta}_f} D_{\mathsf{KL}}(\hat{r}(\mb{\xi})\|p_{\mathsf{KRnet}}(\mb{\xi};\mb{\theta}_f)),
\end{equation*}
where $h(\mb{\xi})$ is defined the same way as in the previous section on a compact support $B_p$ that is slightly larger than $\Omega_p$. 
Again, minimizing the KL divergence is equivalent to minimizing the cross entropy between $\hat{r}$ and $p_{\mathsf{KRnet}}$:
\begin{equation*}
H(\hat{r},p_{\mathsf{KRnet}})=- \int_{B_p} \hat{r}(\mb{\xi}) \log p_{\mathsf{KRnet}}(\mb{\xi};\mb{\theta}_f) d\mb{\xi}.
\end{equation*}
and we approximate the cross entropy using the importance sampling technique:
\begin{equation}\label{eqn:ce_approx}	H(\hat{r},p_{\mathsf{KRnet}}) \approx -\frac{1}{m}\sum_{i=1}^{m} \frac{\hat{r}(\mb{\xi}^{(i)})}{p_{\mathsf{KRnet}}(\mb{\xi}^{(i)};\hat{\mb{\theta}}_f)}\log p_{\mathsf{KRnet}}(\mb{\xi}^{(i)};\mb{\theta}_f),
\end{equation}
where $\hat{\mb{\theta}}_f$ is specified in Algorithm \ref{alg_das_mg}.

\subsection{Algorithm}
We now use the DAS-G strategy presented in \cite{tang_das} to illustrate the procedure of $\mathrm{DAS}^2$, which can also be defined similarly with respect to the DAS-R strategy. 
Given an initial set of collocation points $\mathsf{S}_{\Omega,0}$, the empirical loss defined in \eqref{eq_discrete_residual} is minimized to yield $u_{\mb{\theta}_N^{*,(1)}}$. For $\mb{\theta}_N^{*,(1)}$, one can seek $p_{\mathsf{KRnet}}(\mb{x},\mb{\xi};\mb{\theta}_f^{*,(1)})$ 
by minimizing the cross entropy (see \eqref{eq_discrete_ce}). In this step, uniform samples are used to compute the cross entropy. After the PDF approximation step is finished, a new set of collocation points $\mathsf{S}^g_{\Omega, 1}$ is generated by $p_{\mathsf{KRnet}}(\mb{x},\mb{\xi};\mb{\theta}_f^{*,(1)})$. The training set is refined as $\mathsf{S}_{\Omega, 1} = \mathsf{S}_{\Omega,0} \cup \mathsf{S}^g_{\Omega, 1}$. We then continue to update $u_{\mb{\theta}}$ using $\mb{\theta}_N^{*,(1)}$ as the initial parameters and $\mathsf{S}_{\Omega, 1}$ as the training set, resulting in a refined model. In general, at the $k$-th stage, we minimize the empirical loss on $\mathsf{S}_{\Omega, k-1}$ to get the approximate solution $u_{\mb{\theta}_N^{*,(k)}}$. For PDF approximation, 
we let $p_{\mathsf{KRnet}}(\mb{x},\mb{\xi};\hat{\mb{\theta}}_f) = p_{\mathsf{KRnet}}(\mb{x},\mb{\xi};\mb{\theta}_f^{*,(k-1)})$ for importance sampling in equation \eqref{eq_discrete_ce}. 
Once the PDF model is trained, the training set is refined as $\mathsf{S}_{\Omega, k+1} = \mathsf{S}_{\Omega, k} \cup \mathsf{S}^g_{\Omega, k+1}$. We repeat the procedure to obtain an adaptive algorithm for the refinement of the training set by sampling a joint PDF.

\begin{algorithm}[!htb]
	\caption{$\mathrm{DAS}^2$ based on the joint PDF}
	\label{alg_das_jo}
	\begin{algorithmic}[1]
		\REQUIRE Initial  $p_{\mathsf{KRnet}}(\mb{x},\mb{\xi};\mb{\theta}_f^{(0)})$ , $u_{\mb{\theta}_N^{(0)}}(\mb{x},\mb{\xi})$, maximum epoch number $N_e$, batch size $m$, initial training set $\mathsf{S}_{\Omega, 0} = \{\mb{x}_{0}^{(i)}, \mb{\xi}_{0}^{(i)} \}_{i=1}^{n_r}$.
		\FOR{$k = 0:N_{\rm adaptive}-1$}
		\STATE // Train surrogate models
		\FOR {$i = 1:N_e$}
		\FOR {$j$ steps}
		\STATE Sample $m$ samples from $\mathsf{S}_{\Omega, k}$.
		\STATE Update $u_{\mb{\theta}}(\mb{x},\mb{\xi})$ by descending the stochastic gradient of $J_N(u_{\mb{\theta}})$ (see equation \eqref{eq_discrete_residual}).
		\ENDFOR
		\ENDFOR
		\STATE // Update KRnet
		\FOR {$i = 1:N_e$}
		\FOR {$j$ steps}
		\STATE Sample $m$ samples from $p_{\mathsf{KRnet}}(\mb{x},\mb{\xi};\mb{\theta}_f^{*,(k-1)})$.
		\STATE Update $p_{\mathsf{KRnet}}(\mb{x},\mb{\xi};\mb{\theta}_f)$ by descending the stochastic gradient of $H(\hat{r},p_{\mathsf{KRnet}})$ (see equation \eqref{eq_discrete_ce}).
		\ENDFOR
		\ENDFOR
		\STATE // Refine training set
		\STATE Generate  $\mathsf{S}^g_{\Omega, k+1} \subset \Omega$ with size $n_r$ through $p_{\mathsf{KRnet}}(\mb{x},\mb{\xi};\mb{\theta}_f^{*,(k+1)})$.
		\STATE $\mathsf{S}_{\Omega, k+1} = \mathsf{S}_{\Omega, k} \cup \mathsf{S}^g_{\Omega, k+1}$.
		\ENDFOR	
		\ENSURE $u_{\mb{\theta}_N^*}(\mb{x},\mb{\xi})$
	\end{algorithmic}
\end{algorithm}

\begin{algorithm}[!htb]
	\caption{$\mathrm{DAS}^2$ based on the marginal PDF}
	\label{alg_das_mg}
	\begin{algorithmic}[1]
		\REQUIRE Initial  $p_{\mathsf{KRnet}}(\mb{\xi};\mb{\theta}_f^{(0)})$ , $u_{\mb{\theta}_N^{(0)}}(\mb{x},\mb{\xi})$, maximum epoch number $N_e$, batch size $m$, initial training set $\mathsf{S}_{\Omega_p , 0} = \{ \mb{\xi}_{0}^{(i)} \}_{i=1}^{n_r}$, $m_{\mb{x}}$ samples from $\Omega_s$.
		\FOR{$k = 0:N_{\rm adaptive}-1$}
		\STATE // Train surrogate models
		\FOR {$i = 1:N_e$}
		\FOR {$j$ steps}
		\STATE Sample $m$ samples from $\mathsf{S}_{\Omega_p, k}$.
		\STATE Update $u_{\mb{\theta}}(\mb{x},\mb{\xi})$ by descending the stochastic gradient of $J_N(u_{\mb{\theta}})$ (see equation \eqref{eq_discrete_residual_rd}).
		\ENDFOR
		\ENDFOR
		\STATE // Update KRnet
		\FOR {$i = 1:N_e$}
		\FOR {$j$ steps}
		\STATE Sample $m$ samples from $p_{\mathsf{KRnet}}(\mb{\xi};\mb{\theta}_f^{*,(k-1)})$.
		\STATE Update $p_{\mathsf{KRnet}}(\mb{\xi};\mb{\theta}_f)$ by descending the stochastic gradient of $H(\hat{r},p_{\mathsf{KRnet}})$ (see equation \eqref{eqn:ce_approx}).
		\ENDFOR
		\ENDFOR
		\STATE // Refine training set
		\STATE Generate  $\mathsf{S}^g_{\Omega_p, k+1} \subset \Omega_p$ with size $n_r$ through $p_{\mathsf{KRnet}}(\mb{\xi};\mb{\theta}_f^{*,(k+1)})$.
		\STATE $\mathsf{S}_{\Omega_p, k+1} = \mathsf{S}_{\Omega_p, k} \cup \mathsf{S}^g_{\Omega_p, k+1}$.
		\ENDFOR	
		\ENSURE $u_{\mb{\theta}_N^*}(\mb{x},\mb{\xi})$
	\end{algorithmic}
\end{algorithm}

For simplicity and clarity, we focus on the adaptivity of $\mathsf{S}_{\Omega}$ and the treatment of the boundary points can be found in \cite{tang_das} (section 4.3).  The deep adaptive sampling algorithm for surrogate modeling is summarized in Algorithm \ref{alg_das_jo}, where $N_{\text{adaptive}}$ is a
given number of maximum adaptivity iterations, $m$ is the batch size for stochastic gradient, and $N_e$ is the number of epochs for training $u_{\mb{\theta}}(\mb{x},\mb{\xi})$ and $p_{\mathsf{KRnet}}(\mb{x},\mb{\xi};\mb{\theta}_f)$. The algorithms consist of three steps in one loop: training surrogate models, updating KRnet and refining the training set. 
The same procedure can be applied to the marginal PDF, which results in Algorithm \ref{alg_das_mg}.

\section{Analysis}
Inspired by the literature \cite{de2022error, de2022generic}, we include some preliminary analysis of $\mathrm{DAS}^2$. 
We first establish the relationship between the loss functional and its discretization at the optimal model parameters for a certain training set. 
For the ideal case, we show that the expectation of the discretized loss functional does not increase at the optimal model parameters given by two adjacent adaptivity iterations. Before presenting the analysis, the following assumptions are introduced.

\begin{assump}
	[\cite{de2022error}]
	\label{assump_lip_operator}
	Let $\mb{\theta} \in \Theta = [-a,a]^D$ be the trainable parameters of $u_{\mb{\theta}}$ where $a > 0$ is a constant. Assume that two operators $\mathscr{M}_1: \mb{\theta} \mapsto J_{r,N}$ and $\mathscr{M}_2: \mb{\theta} \mapsto J_{r}$ are Lipschitz continuous in the $\ell_{\infty}$ sense with Lipschitz constant $\mathfrak{L}$ for $\mb{\theta} \in \Theta$. 
\end{assump}

\begin{assump}
	[\cite{de2022error}]
	\label{assump_bounded}
	Let $c > 0$ be a constant that is independent of $\Theta$. Assume that $J_{r,N}\in [0,c]$ for all $\mb{\theta} \in \Theta$. 
\end{assump}

\begin{assump}[\cite{tang_das}]\label{assump_pdf}
	Assume that $p_{\mathsf{KRnet}}(\bx,\mb{\xi}; \btheta_f^{*,(k)})$ is the optimal candidate for the change of measure for problem \eqref{empminimizer} at $k$-th stage 
	\begin{equation*}
		p_{\mathsf{KRnet}}(\bx,\mb{\xi}; \btheta_f^{*,(k)}) = c_k r^2(\bx,\mb{\xi}; \btheta_N^{*,(k)}),
	\end{equation*}
	where $\btheta_N^{*,(k)}$ is the minimizer in \eqref{empminimizer} and $\btheta_f^{*,(k)}$ is the minimizer in \eqref{eq_kl_krres} given $\btheta_N^{*,(k)}$, and $$c_k=1/\int_{\Omega}r^2(\bx,\mb{\xi}; \btheta_N^{*,(k)})d\bx d\mb{\xi}$$ is the normalization constant. 
\end{assump}

If the collocation points are independently and identically distributed according to a given probability distribution, then $J_r(u_{\mb{\theta}_N^*})$ can be bounded by the discrete residual with high probability, which is stated as follows. 
\begin{thm}\label{thm_bounded_by_residual}
	Suppose that Assumption \ref{assump_lip_operator} and Assumption \ref{assump_bounded} are satisfied and the boundary loss is zero. Let $\mb{\theta}_N^*$ be a minimizer of $J_{r,N}$ where the collocation points are independently drawn from a given probability distribution. Given $\varepsilon \in (0,1)$, the following inequality holds
	\begin{equation*}
		J_r(u_{\mb{\theta}_N^{*}}) \leq \varepsilon^2 + J_{r,N}(u_{\mb{\theta}_N^{*}})
	\end{equation*}
	with probability at least $1 - (4a \mathfrak{L}/\varepsilon^2)^D \mathrm{exp}(-N_r \varepsilon^4/2c^2 )$.
\end{thm}

The expectation of the discrete residual at two adjacent adaptivity stages satisfies the following property.
\begin{thm}\label{thm_error_beha}
	Under the same conditions of Theorem \ref{thm_bounded_by_residual}, suppose that Assumption \ref{assump_pdf} is satisfied. Assume that 
	\begin{equation*}
		J_{r,N}(u_{\mb{\theta}_N^{*,(k)}}) = \frac{1}{N_r} \sum\limits_{i=}^{N_r}  \frac{r^2(\mb{x}^{(i)},\mb{\xi}^{(i)};\mb{\theta}_N^{*,(k)})}{p_{\mathsf{KRnet}}(\mb{x}^{(i)},\mb{\xi}^{(i)}; \btheta_f^{*,(k-1)})}, 
	\end{equation*}
    where each $(\mb{x}^{(i)}, \mb{\xi}^{(i)})$ is drawn from $p_{\mathsf{KRnet}}(\mb{x}^{(i)},\mb{\xi}^{(i)}; \btheta_f^{*,(k-1)})$, 
	then the following inequality holds
	\begin{equation*}
		\mathbb{E}(J_{r,N}(u_{\mb{\theta}_N^{*,(k+1)}})) \leq \mathbb{E}(J_{r,N}(u_{\mb{\theta}_N^{*,(k)}})).
	\end{equation*}
\end{thm}

The proofs of Theorem \ref{thm_bounded_by_residual} and Theorem \ref{thm_error_beha} can be found in the Appendices. 
Theorem \ref{thm_bounded_by_residual} provides the relationship between the residual and the discrete residual, which is similar to the results in \cite{de2022error, de2022generic}. The analysis in Theorem \ref{thm_bounded_by_residual} is not restricted to the linear differential equations, while the results in \cite{tang_das} only involve non-parametric linear differential equations. From the above analysis, if the number of parameters in $u_{\mb{\theta}}$ and the number of samples are properly chosen, then the residual is bounded by the discrete residual (i.e., loss) with high probability. 
In Theorem \ref{thm_error_beha}, we consider the DAS-R strategy for simplicity. Although the analysis of the residual behavior during the adaptive procedure is restricted to DAS-R, it can provide a perspective to understand the mechanism of $\mathrm{DAS}^2$. We note that quantifying the decay of the error is not straightforward since it depends on the optimization procedure of deep neural networks. However, obtaining the convergence rate of such an optimization problem is still an open question.

\section{Numerical study}\label{sec_numexp}
In this section, we conduct four numerical experiments (including two parametric ordinary differential equations and two parametric partial differential equations) to demonstrate the effectiveness of the proposed method, where different types of problems under different parametric settings are studied. Two types of neural network structures are considered for the surrogate model $u_{\mb{\theta}}(\mb{x}, \mb{\xi})$: one is the feedforward neural network with inputs $\mb{x}$ and $\mb{\xi}$, and the other one is the structure given by DeepONet \cite{lu2021learning} where the problem is treated as an operator learning problem. The choice of the sampling strategy depends on the problem properties instead of the model structure. 
For comparison, we also test the performance of some baseline sampling strategies, such as the residual-based adaptive refinement (RAR) method \cite{lu2021deepxde,wu2023comprehensive} and the quasi-random sampling (QRS) method implemented in the SciPy module \cite{virtanen2020scipy}. 
We use DAS-G for all numerical experiments since it is more robust than DAS-R \cite{tang_das}.
The code of this study will be released on GitHub once the paper is accepted.

\subsection{A one-dimensional parametric ordinary differential equation}
We start with the following one-dimensional parametric ordinary differential equation (ODE)
\begin{equation*}
	\frac{\text{d}u}{\text{d}x}=\xi u,\quad u(0,\xi)=u_0,\quad x\in[0,1],
\end{equation*}
where $\xi \in \Omega_p =  [-3,3]$, i.e., $C_1 = -3$ and $C_2 = 3$ (see section \ref{sec_probset}), and the initial condition is set to $u_0 = 1$. The exact solution is 
\begin{equation}\label{eq:ex1_true}
	u(x,\xi)=u_0 e^{\xi x}.
\end{equation}
This is a widely used test problem for polynomial chaos methods in uncertainty quantification \cite{xiu2002wiener}.

We use a six-layer fully connected neural network to construct a surrogate model $u_{\mb{\theta}}(x,\xi)$ as the approximation solution of the parametric ODE, where each hidden layer has 32 neurons. 
For KRnet, we set $K = 2$ and take $L=6$ affine coupling layers. For each affine coupling layer, a two-layer fully connected neural network is used, where each hidden layer has 24 neurons. The maximum epoch number for training both $u_{\mb{\theta}}(x,\xi)$ and $p_{\mathsf{KRnet}}(x,\xi;\mb{\theta}_f)$ is set to $N_e=3000$. In this test problem, the ADAM optimizer \cite{kingma2017adam} is employed for all training processes. The learning rate for the ADAM optimizer is set to 0.0001, and the batch size is set to $m=1000$. For $\mathrm{DAS}^2$, we use the joint PDF for sampling. The collocation points in the initial training set are uniform samples, $n_r=1000$ is set during the adaptive sampling procedure, and the number of adaptivity iterations is set to $N_{\rm adaptive}=6$. For the uniform sampling strategy, the maximum epoch number is set to be the same as the total number of epochs of $\mathrm{DAS}^2$, and the number of samples is set to $|\mathsf{S}_{\Omega}|=6000$ (the same as that of $\mathrm{DAS}^2$). To assess the effectiveness of our $\mathrm{DAS}^2$ method, we generate a uniform meshgrid with size $256\times 256$ in the spatial-parametric space $[0,1]\times [-3,3]$ and compute the mean square error on these grid points.

\begin{figure}[!htb]
	\centering
	\subfloat[][The error with respect to epoch.]{\includegraphics[width=.45\textwidth]{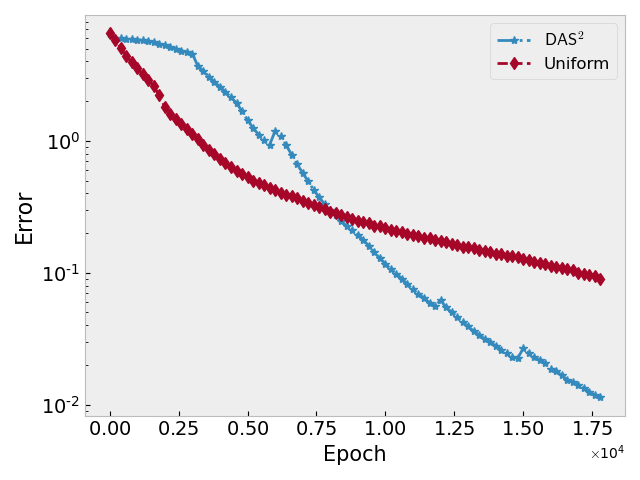}}
	\subfloat[][The error of $\mathrm{DAS}^2$ at certain adaptivity iteration steps.]{\includegraphics[width=.45\textwidth]{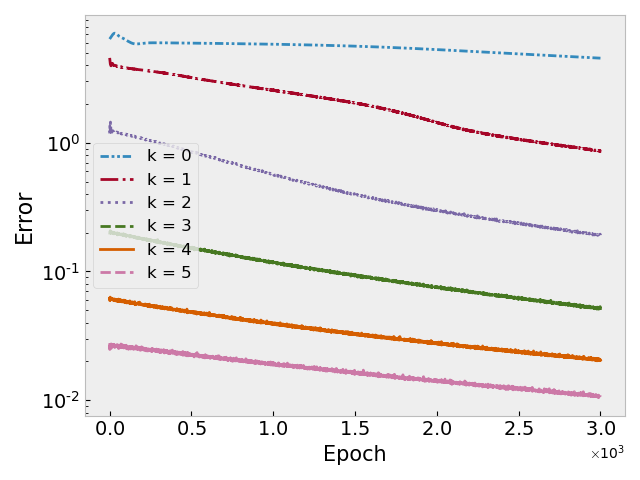}}
	\caption{The errors for the parametric ODE test problem.}
	\label{fig:ode_err}
\end{figure}

\begin{figure}[!htb]
	\centering
	\subfloat[][The evolution of $\mathsf{S}_{\Omega,k}^g$ in $\mathrm{DAS}^2$.]{\includegraphics[width=.63\textwidth]{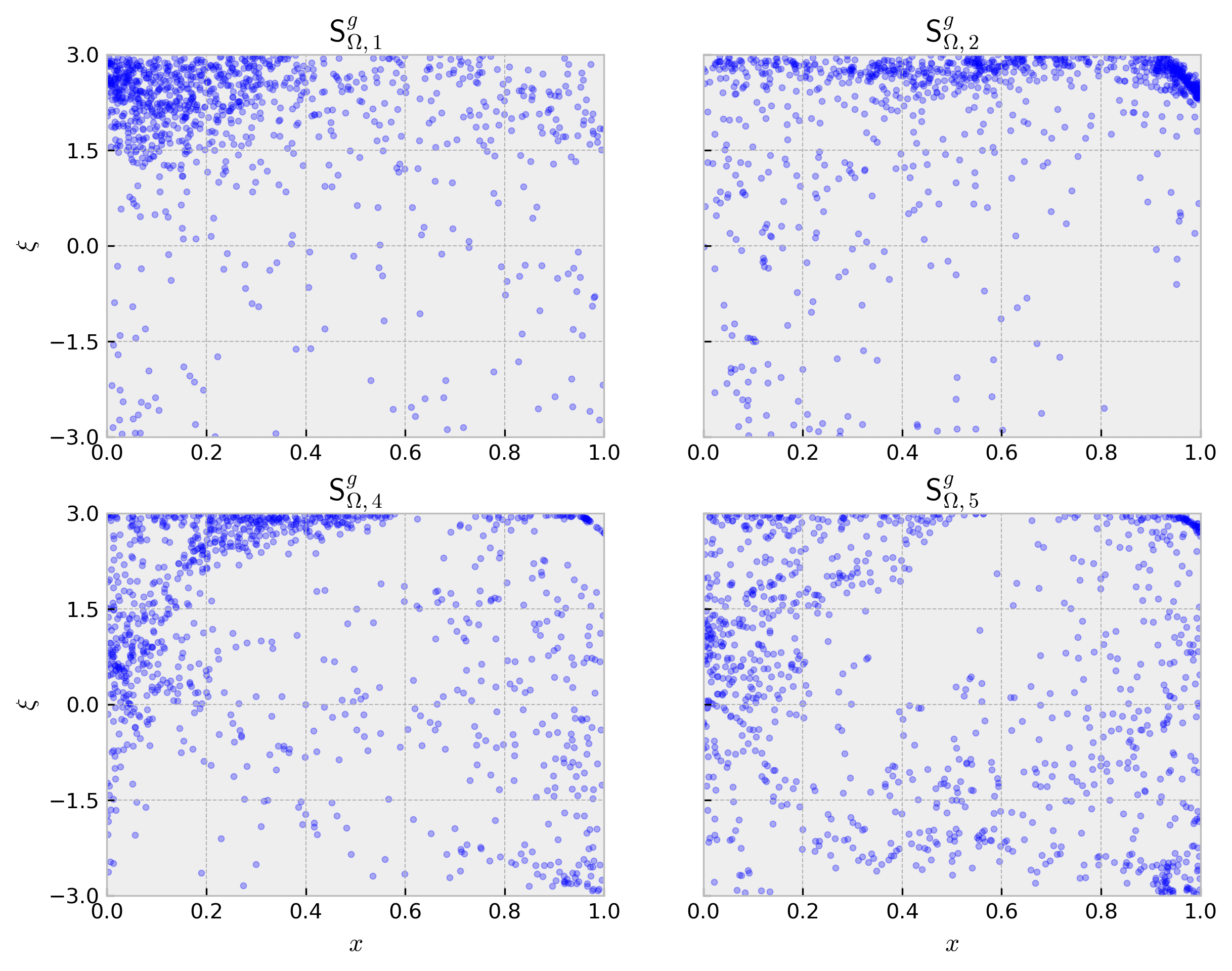}}
	\subfloat[][$u(x,\xi)$.]{\includegraphics[width=.32\textwidth]{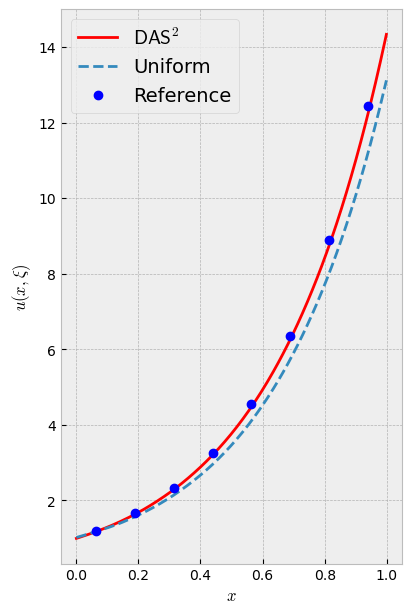}}
	\caption{The results of the parametric ODE test problem. Left: The evolution of $\mathsf{S}_{\Omega,k}^g$ in $\mathrm{DAS}^2$; Right: The exact solution and the approximate solutions with different sampling strategies for $\xi = 2.725$.}
	\label{fig:ex1_sampling_pre}
\end{figure}

In Figure~\ref{fig:ode_err}, we plot the approximation error given by different sampling strategies with respect to epoch in the left plot and the error evolution of $\mathrm{DAS}^2$ at different adaptivity iteration steps in the right plot. In terms of the number of epochs, the error of $\mathrm{DAS}^2$ decays more quickly than the uniform sampling method. 
The approximation error of $\mathrm{DAS}^2$ drops as the adaptivity iteration step $k$ increases. It can be found from equation \eqref{eq:ex1_true} that $u(x,\xi)$ grows exponentially with respect to $x$ and $\xi$. When $\xi$ is near $3$, the solution increases dramatically. To capture this information, more samples are located in the area that $\xi$ is near $3$. Figure \ref{fig:ex1_sampling_pre}(a) shows the evolution of $\mathsf{S}_{\Omega}^g$ of $\mathrm{DAS}^2$ with respect to adaptivity iterations $k = 1,2,4,5$ ($|\mathsf{S}_{\Omega,k}^g| = 1000$), where $\mathsf{S}_{\Omega,1}^g$ indicates that large point-wise residuals are located in the upper left corner of the $x$-$\xi$ plane. After the set of collocation points is augmented by $\mathsf{S}_{\Omega,1}^g$, the residual profile changes as shown in $\mathsf{S}_{\Omega,2}^g$. Such a pattern is repeated until $\mathsf{S}_{\Omega,k}^g$ is near a set of uniform samples. Figure \ref{fig:ex1_sampling_pre}(b) shows the exact solution, the solution obtained by $\mathrm{DAS}^2$ and the solution obtained by the uniform sampligng strategy for $\xi = 2.725$. It is seen that $\mathrm{DAS}^2$ yields a more accurate approximation.


\subsection{Operator learning for a dynamical system with high-dimensional parameters}
Next we consider the following dynamical system
\begin{equation}\label{eq:op_1d}
	\left\{
	\begin{array}{r l}
		\dfrac{\text{d}u(x,\mb{\xi})}{\text{d}x}= e^{-D\|\bm{\xi}-\mb{0.5}\|^2}f(x,\mb{\xi}),& x\in[0,1],\\
		u(0,\mb{\xi})=0, & 
	\end{array}
	\right.
\end{equation}
where $D$ is a fixed parameter, and $\mb{\xi} \in \Omega_p = [-M, M]^d$, i.e., $C_1 = -M$ and $C_2 = M$. The goal is to learn the solution operator from $f$ to the solution $u$ without any paired input-output data when $f$ is sampled from a given function space. This example without $\mathrm{exp}(-D\|\bm{\xi}-\mb{0.5}\|^2)$ is used to test the performance of DeepONet \cite{lu2021learning, wang2021learning}. Here, we add a term $\mathrm{exp}(-D\|\bm{\xi}-\mb{0.5}\|^2)$ to the right-hand side to make this problem more challenging.  We assume that $f$ is drawn from the space spanned by orthogonal (e.g. Chebyshev) polynomials as studied in \cite{lu2021learning}. Let $T_i$ be Chebyshev polynomials of the first kind. We define the orthogonal polynomials of degree $d$ as:
\begin{equation*}
	V_{\text{poly}}=\left\{ \sum_{i=0}^{d-1} \xi_i T_i(x):\left|\xi_i\right|\leq M   \right\}.
\end{equation*}
This function space is parameterized with $\mb{\xi}= [\xi_0,\xi_1,...,\xi_{d-1}]^{\mathsf{T}}$. 
Given a realization of $\mb{\xi}$, we can generate a continuous function $f$ as the following form
\begin{equation*}
	f(x,\mb{\xi})=\sum_{i=0}^{d-1} \xi_i T_i(x).
\end{equation*}
In this example, the parametric solution $u(x,\mb{\xi})$ is approximated by 
\begin{equation}\label{eq:deeponet}
	u_{\mb{\theta}}(x,\mb{\xi}) \approx \sum\limits_{i=1}^l q_{\mb{\theta}_1}^{(i)}(x) t_{\mb{\theta}_2}^{(i)}(\mb{\xi}) + b_0,
\end{equation}
where $q_{\mb{\theta}_1}^{(i)}$ and $t_{\mb{\theta}_2}^{(i)}$ are $i$-th outputs of two neural networks $q$ (parameterized with $\mb{\theta}_1$) and $t$ (parameterized with $\mb{\theta}_2$) respectively, both of which have $l$ outputs, and $b_0 \in \xs{R}$ is a bias to be trained. Denoting the whole parameters in \eqref{eq:deeponet} by $\mb{\theta} = \{\mb{\theta}_1, \mb{\theta}_2, b_0 \}$ for short. 

The experimental setup is as follows. We set $M=1,d=8,D=6$. $q_{\mb{\theta}_1}(x)$ and $t_{\mb{\theta}_2}(\mb{\xi})$ are both five-layer fully connected neural networks and each hidden layer has 50 neurons. For $\mathrm{DAS}^2$, we use a marginal PDF for adaptive sampling because the singularity is mainly in the parametric space. To compute the marginal PDF and the loss functional, we use $m_{\mb{x}}=100$ uniform grid points in $[0,1]$ to discretize the integral in equation \eqref{eq_mgpdf_discrete} and \eqref{eq_discrete_residual_rd}, {in other words, adaptive sampling is not considered in the physical space since the low regularity is from the parametric space. 
We set $K=4$ and the configuration for the affine coupling layer is the same as the previous example. 
The number of epochs for training both $u_{\mb{\theta}}(x,\mb{\xi})$ and $p_{\mathsf{KRnet}}(\mb{\xi};\mb{\theta}_f)$ is set to $N_e=3000$. The learning rate for the ADAM optimizer is set to 0.0001, and the batch size is set to $m = 5000$. The numbers of adaptivity iterations is set to $N_{\rm adaptive}=5$. For the uniform sampling strategy, we generate $f(x,\mb{\xi})$ with each $\xi_i \sim \mathsf{Uni}(-M,M)$ where $\mathsf{Uni}(-M,M)$ is the uniform distribution on $[-M, M]$. To measure the quality of approximation, we generate a validation set, which contains $10000$ uniformly distributed points in $[-M,M]^d$ and $10000$ points in the $d$-dimensional ball centered at $\mb{0.5}$ with radius $0.5$. To compute the reference solution, we employ the classical Runge-Kutta45 (RK45) method to solve the ODE for each function $f$ with a certain $\mb{\xi}$.

\begin{figure}[!htb]
    \centering
    \subfloat[][The error w.r.t. sample size $\vert \mathsf{S}_{\Omega_p} \vert$.]{\includegraphics[width=.45\textwidth]{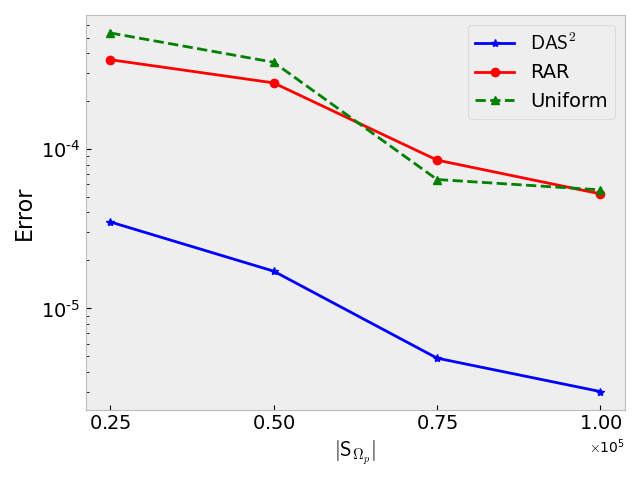}}
    \subfloat[][The error evolution with $|\mathsf{S}_{\Omega_p}|=1\times 10^5$.]{\includegraphics[width=.45\textwidth]{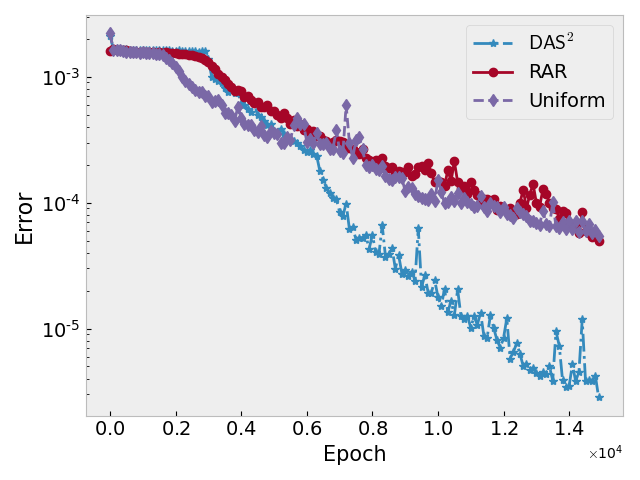}}
    \caption{Approximation errors for the operator learning problem.}
    \label{fig:ex3_comparison}
\end{figure}

\begin{figure}[!htb]
    \centering
    \includegraphics[width=0.55\textwidth]{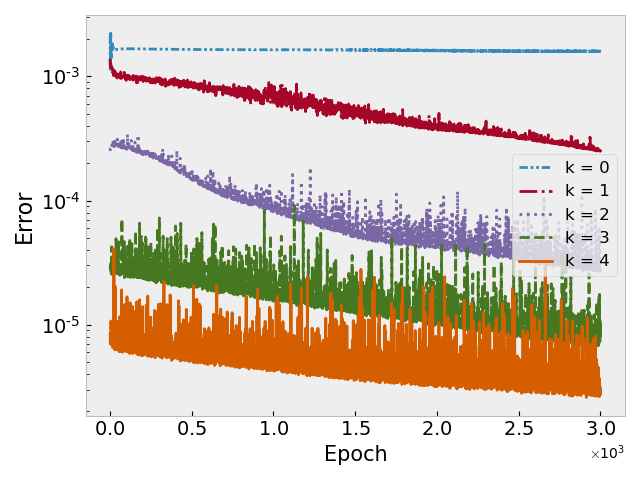}
    \caption{The error evolution of $\mathrm{DAS}^2$ at different adaptivity iteration steps for the operator learning problem. $|\mathsf{S}_{\Omega_p}|=1\times 10^5.$}
    \label{fig:ex3_das_stage}
\end{figure}

\begin{table}[!htb]
    \caption{The operator learning problem: inference time and error for different $\vert \mathsf{S}_{\Omega_p} \vert$ and sampling strategies. The computing time of RK45 is 105 seconds.}
    \centering	
    \begin{small}
        \begin{tabular}{cccccccccc}  
            \toprule 
            \multicolumn{2}{c}{\multirow{2}*{\diagbox{sampling strategy}{$\vert \mathsf{S}_{\Omega_p} \vert$}}} 
            & \multicolumn{1}{c}{$2.5\times 10^4$}& & \multicolumn{1}{c}{$5\times 10^4$} & & \multicolumn{1}{c}{$7.5\times 10^4$}  & & \multicolumn{1}{c}{$1\times 10^5$} \\ 
            \\
            \hline
            \multicolumn{2}{l}{Uniform (0.006s)} &   $5.4\times 10^{-4}$ & &   $3.5\times 10^{-4}$ & &   $6.4\times 10^{-5}$ & & $5.5\times 10^{-5}$\\
            \multicolumn{2}{l}{RAR (0.006s)} &    $3.6\times 10^{-4}$ & &  $2.6\times 10^{-4}$ & &   $8.5\times 10^{-5}$ & & $5.2\times 10^{-5}$ \\
            \multicolumn{2}{l}{$\mathrm{DAS}^2$ (0.03s)}    & $3.5\times 10^{-5}$ & & $1.7\times 10^{-5}$ &  & $4.9\times 10^{-6}$ & & $3.0\times 10^{-6}$ \\
            \bottomrule
        \end{tabular}
    \end{small}
    \label{oplearning_res}
\end{table}

In Figure~\ref{fig:ex3_comparison}, we plot the mean square error of different sampling strategies with respect to the sample size $\vert \mathsf{S}_{\Omega_p} \vert$ in the left plot and with respect to the number of epochs in the right plot. For $\mathrm{DAS}^2$, the numbers of collocation points in $\mathsf{S}_{\Omega_p,k}^g$ ($k = 1,2,3,4$) are set to $n_r = 5 \times 10^3, 1 \times 10^4, 1.5 \times 10^4, 2 \times 10^4$ for $\vert \mathsf{S}_{\Omega_p} \vert = 2.5 \times 10^4, 5 \times 10^4, 7.5 \times 10^4, 1 \times 10^5$ respectively. For the uniform sampling strategy, the model is trained with $1.5 \times 10^4$ epochs to match the total number of epochs of $\mathrm{DAS}^2$. For the heuristic method RAR, the numbers of collocation
points in $\mathsf{S}_{\Omega_p,k}^g$ ($k = 1,2,3,4$) are set to $n_r = 2.5 \times 10^3, 5 \times 10^3, 7.5 \times 10^3, 1 \times 10^4$ for $\vert \mathsf{S}_{\Omega_p} \vert = 2.5 \times 10^4, 5 \times 10^4, 7.5 \times 10^4, 1 \times 10^5$  respectively. From the left plot of Figure~\ref{fig:ex3_comparison}, it can be seen that $\mathrm{DAS}^2$ improves the accuracy significantly compared to the uniform sampling strategy and RAR.
The right plot of Figure~\ref{fig:ex3_comparison} shows that as the number of epochs increases, especially from the start of the third adaptivity iteration, the error of $\mathrm{DAS}^2$ decreases much faster than those of 
uniform sampling and RAR. 
Figure~\ref{fig:ex3_das_stage} shows the errors of $\mathrm{DAS}^2$ at each adaptivity iteration step $k$. It is seen that the error drops dramatically after we refine the solution using $\mathsf{S}_{\Omega_p,1}^g$ and $\mathsf{S}_{\Omega_p,2}^g$. Table \ref{oplearning_res} shows the inference time and the errors for the uniform sampling strategy, RAR and $\mathrm{DAS}^2$. As a surrogate model, the inference time of $\mathrm{DAS}^2$ is much less than that of RK45, which is desired. It can be seen that the inference time of $\mathrm{DAS}^2$ is more than that of RAR. However, the errors of the uniform sampling strategy and RAR are much larger than that of $\mathrm{DAS}^2$ since the uniform sampling strategy and RAR are not able to accurately discretize the loss functional for this low-regularity high-dimensional problem \cite{wu2023comprehensive}. From Table \ref{oplearning_res}, it is clear that $\mathrm{DAS}^2$ is one order of magnitude more accurate than RAR and the uniform sampling strategy.

\begin{figure}[!htb]
    \centering
    \includegraphics[width=0.85\textwidth]{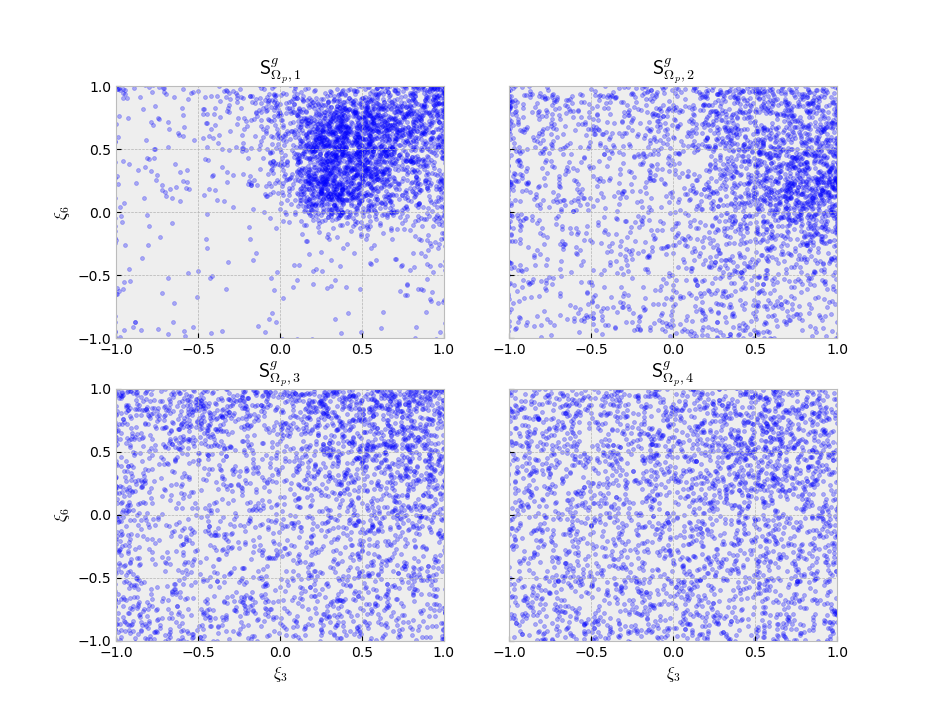}
    \caption{The evolution of $\mathsf{S}_{\Omega_p,k}^g$ in $\mathrm{DAS}^2$ for the operator learning problem, $|\mathsf{S}_{\Omega_p}|=1\times 10^5$.}
    \label{fig:ex3_sampling}
\end{figure}

\begin{figure}[!htb]
    \centering
    \includegraphics[width=1.0\textwidth]{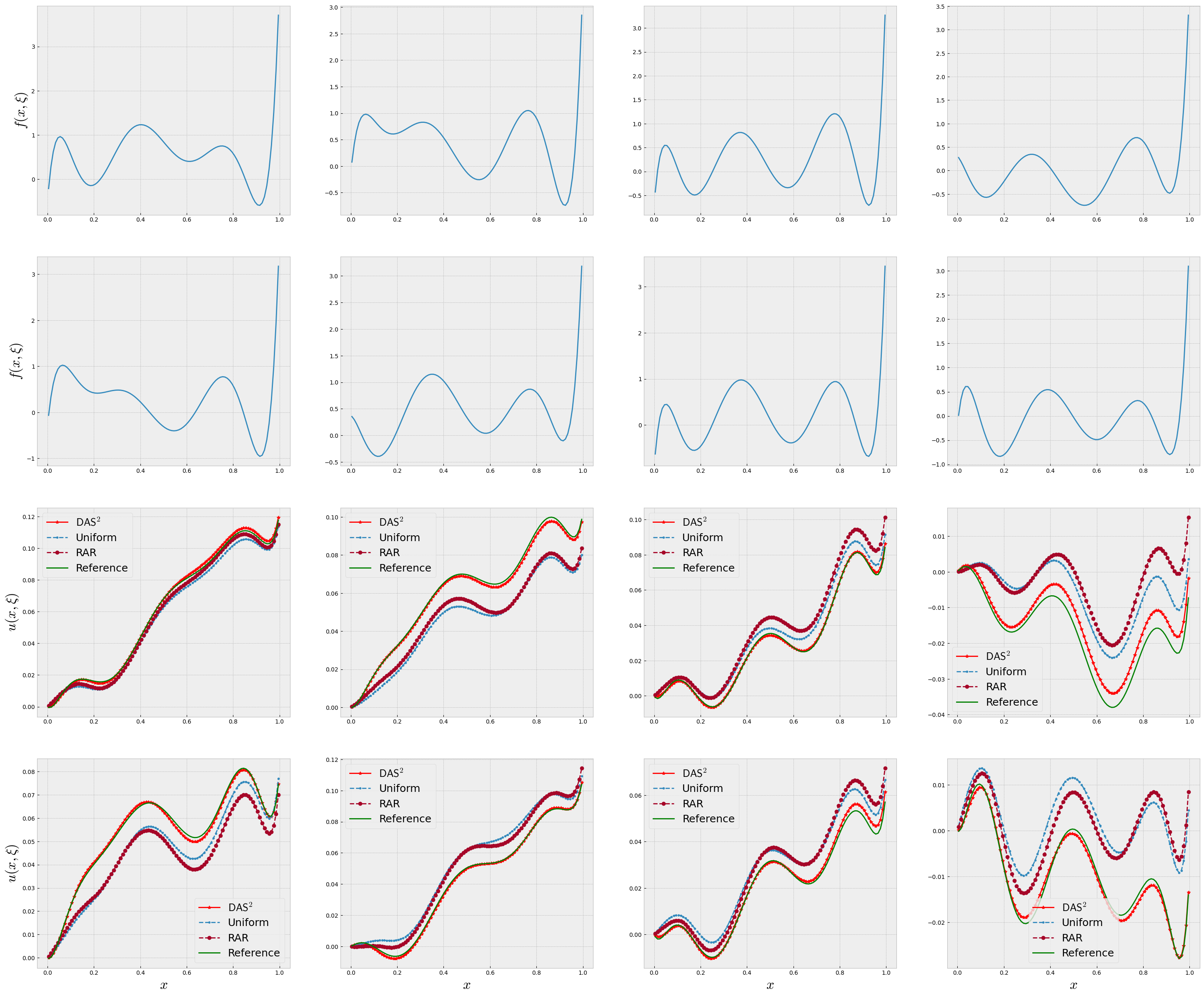}
    \caption{Solutions of the operator learning problem: the first two rows above show $f(x,\mb{\xi})$ with different realizations of $\mb{\xi}$, the two rows below show the corresponding solutions $u(x,\mb{\xi})$.}
    \label{fig:ex3_presentation}
\end{figure}

Figure~\ref{fig:ex3_sampling} shows $3000$ samples from $\mathrm{DAS}^2$ for the four adaptivity iterations, where the components $\xi_3$ and $\xi_6$ are used for visualization. We have also checked the other components, and no significantly different results were found. $\mathsf{S}_{\Omega_p ,1}^g$ shows that the error profile has a peak around $\hat{\mb{\xi}}=\mb{0.5}$ which matches the fact that there is a decay term with respect to $\mb{\xi}$ in equation \eqref{eq:op_1d}. After the training set is augmented with $\mathsf{S}_{\Omega_p,1}^g$, the error profile becomes more flat as shown by the distribution of $\mathsf{S}_{\Omega_p,2}^g$. This is expected since more training samples are added to the neighborhood of $\mb{0.5}$ where the error should be reduced. Figure~\ref{fig:ex3_presentation} shows $u(x,\mb{\xi})$ corresponding to different $\mb{\xi}$ obtained by $\mathrm{DAS}^2$, RAR and the uniform sampling method. The realizations of $\mb{\xi}$ we choose for visualization are randomly drawn from the $d$-dimensional ball centered at $\mb{0.5}$ with radius 0.5, since $u(x,\mb{\xi})$ is close to zero when $\mb{\xi}$ is far away from $\mb{0.5}$ due to the decay term in problem \eqref{eq:op_1d}. As shown in Figure~\ref{fig:ex3_presentation}, for different $\mb{\xi}$ the solutions $u(x,\mb{\xi})$ obtained by $\mathrm{DAS}^2$ are much more accurate than those given by RAR and uniform sampling. 

\subsection{Surrogate modeling for an optimal control problem with geometrical parametrization}
In this test case, we are going to build a surrogate model for the following parametric optimal control problem: 
\begin{linenomath*}\begin{equation}
        \label{eq:pocp}
        \left\{\begin{aligned} 
            &\min _{y(\mb{x},\mb{\xi}), u(\mb{x},\mb{\xi})} J\left(y\left(\mb{x},\mb{\xi}\right), u\left(\mb{x},\mb{\xi}\right)\right)=\frac{1}{2}\left\|y(\mb{x},\mb{\xi})-y_{d}(\mb{x},\mb{\xi})\right\|_{2,\Omega}^{2}+\frac{\alpha}{2}\left\|u(\mb{x},\mb{\xi})\right\|_{2,\Omega}^{2}, \\
            &\text{subject to}  
            \left\{\begin{aligned}
                -\Delta y(\mb{x},\mb{\xi}) &= u(\mb{x},\mb{\xi}) &&\text {in} \; \Omega, \\ 
                y(\mb{x},\mb{\xi}) &=1 &&\text {on}\; \partial \Omega,\\
            \end{aligned}\right.\\ 
            &\text{and} \quad u_a \leq u(\mb{x},\mb{\xi}) \leq u_b \quad \text {a.e. in}\; \Omega,
        \end{aligned}\right.
\end{equation}\end{linenomath*}
where $\mb{\xi}=(\xi_1, \xi_2)$ represents the geometrical and desired state parameters. The parametric computational domain (also depending on $\mb{\xi}$) is $\Omega = ([0,2]\times[0,1])\backslash \mathbb{B}((1.5,0.5), \xi_1)$ which is illustrated in Figure~\ref{fig:ex4_illustration} and the desired state is given by
\begin{linenomath*}\begin{equation*}
        y_d(\mb{x},\mb{\xi}) = 
        \begin{cases}
            1 & \text {in}\; \Omega_{1}=[0,1]\times[0,1], \\ 
            \xi_2 & \text {in} \; \Omega_{2} =([1,2]\times[0,1])\backslash \mathbb{B}((1.5,0.5),\xi_1),
        \end{cases}
\end{equation*}\end{linenomath*}
where $\mathbb{B}((1.5,0.5), \xi_1)$ is a ball of radius $\xi_1$ with center $(1.5,0.5)$. We set $\alpha=0.001, u_a=0, u_b=10$, and the domain for the parameter to be $\mb{\xi} \in \Omega_p = [0.05,0.45]\times[0.5,2.5]$. 
\begin{figure}[!htb]
    \centering
    \includegraphics[width=0.5\textwidth]{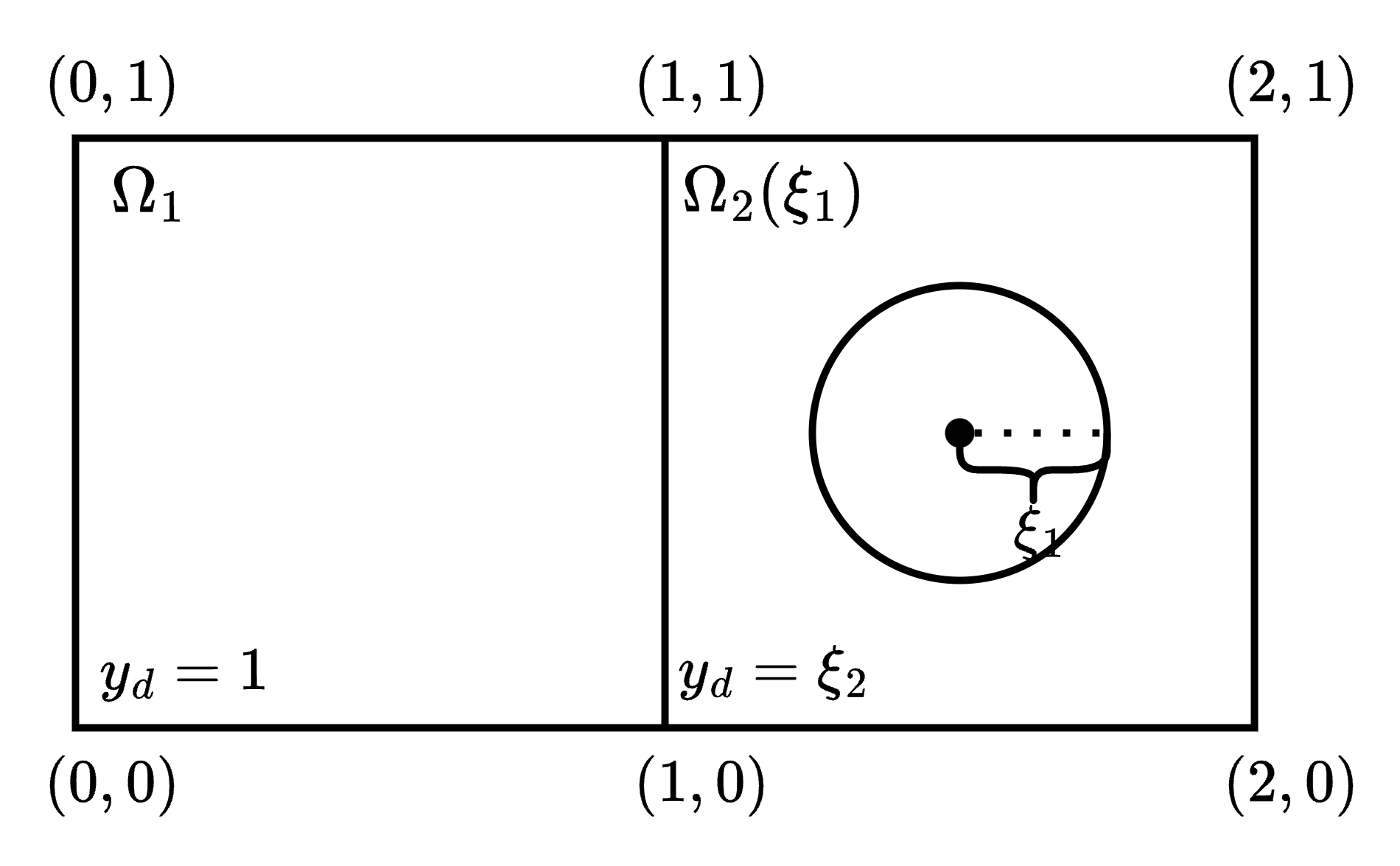}
    \caption{The parametric computational domain $\Omega$.}
    \label{fig:ex4_illustration}
\end{figure}
This test problem is related to the application of local hyperthermia treatment of cancer, which is inspired by the literature \cite{karcher2018certified, negri2013reduced}. The background of this test problem is that we expect to accomplish a specific temperature field in the tumor area and another temperature field in the non-lesion area by heat source control. The circle represents a certain body organ where the tumor area is.  We intend to seek an effective surrogate model of the optimal heat source control for different expected temperature fields and organ shapes (i.e. different $\mb{\xi}$).

As studied in \cite{yin2023aonn}, one can use the necessary conditions for the minimizer of \eqref{eq:pocp} to find the optimal solution to the parametric optimal control problem. That is, we solve the KKT system of \eqref{eq:pocp} to find its minimizer, which is a parametric PDE system as follows
\begin{equation}
    \left\{
    \begin{aligned}
        -\Delta y(\mb{x},\mb{\xi})=u(\mb{x},\mb{\xi}) \quad & \text{in}\;\Omega, \\
        y(\mb{x},\mb{\xi})=1 \quad & \text{on}\; \partial\Omega, \\
        -\Delta p(\mb{x},\mb{\xi})=y(\mb{x},\mb{\xi})-y_d(\mb{x},\mb{\xi})\quad & \text{in}\; \Omega, \\ 
        p(\mb{x},\mb{\xi})=0 \quad & \text{on}\;\partial \Omega, \\
        u(\mb{x},\mb{\xi})=-\frac{1}{\alpha}\text{P}_{[u_a,u_b]}\left(p\left(\mb{x},\mb{\xi}\right)\right) \quad  & \text{in}\; \Omega,\\
    \end{aligned}
    \right.
    \label{eq:ex4_kkt}
\end{equation}
where $p(\mb{x},\mb{\xi})$ is the adjoint variable and
\begin{equation*}
    \text{P}_{[u_a,u_b]}\left(p\left(\mb{x},\mb{\xi}\right)\right) = \left\{
    \begin{aligned}
        u_b, \quad & \text{if}\;u_b < p(\mb{x},\mb{\xi}),\\
        p(\mb{x},\mb{\xi}), \quad & \text{if}\; u_a \leq p(\mb{x},\mb{\xi}) \leq u_b,\\
        u_a,\quad & \text{if}\; p(\mb{x},\mb{\xi})< u_a.
    \end{aligned}
    \right.
\end{equation*}
Define a length factor function as \cite{yin2023aonn}
\begin{equation*}
    l(\mb{x},\mb{\xi})=x_1(2-x_1)x_2(1-x_2)(\xi_1^2-(x_1-1.5)^2-(x_2-0.5)^2).
    \label{eq:ex4_lf}
\end{equation*}
We choose three six-layer fully connected neural networks $u_{\mb{\theta}_u}(\mb{x},\mb{\xi})$, $y_{\mb{\theta}_y}(\mb{x},\mb{\xi})$ and $p_{\mb{\theta}_p}(\mb{x},\mb{\xi})$, where each hidden layer has 25 neurons. We let
\[
u(\mb{x},\mb{\xi})\approx u_{\mb{\theta}_u}(\mb{x},\mb{\xi}),\,  y(\mb{x},\mb{\xi}) \approx l(\mb{x},\mb{\xi}) y_{\mb{\theta}_y}(\mb{x},\mb{\xi}) + 1,\, p(\mb{x},\mb{\xi}) \approx  l(\mb{x},\mb{\xi})p_{\mb{\theta}_p}(\mb{x},\mb{\xi})
\]
The Dirichlet boundary conditions of $y(\mb{x},\mb{\xi})$ and $p(\mb{x},\mb{\xi})$ are naturally satisfied. We then substitute the defined approximators into equation \eqref{eq:ex4_kkt} to minimize the residual. More details about the discretization of problem \eqref{eq:pocp} can be found in \cite{yin2023aonn}. We here focus on the importance of adaptive sampling for surrogate modeling. 


In this example, we use the joint PDF for sampling and the spatial-parametric space is defined as:
\begin{equation*}
    \begin{aligned}
        \Omega :=\{
        (\mb{x},\mb{\xi})|& 0\leq x_1\leq 2, 0\leq x_2\leq 1,0.05\leq \xi_1 \leq 0.45,0.5\leq \xi_2\leq 2.5, \\
        & (x_1-1.5)^2+(x_2-0.5)^2\geq \xi_1^2
        \}. \\
    \end{aligned}
\end{equation*}
To obtain an accurate approximation, the optimizer for training $u_{\mb{\theta}_u}$, $y_{\mb{\theta}_y}$ and $p_{\mb{\theta}_p}$ is set to be the BFGS method \cite{jorge2006numerical}, followed by the setting in \cite{yin2023aonn}. The number of epochs for solving PDEs is set to $N_e=2000$. 
For KRnet, we set $K=2$ 
and the configuration for the affine coupling layers remains the same as the previous experiment. 
KRnet is trained by the ADAM optimizer with a learning rate 0.0001, where the number of epochs is set to $N_e=2000$. For $\mathrm{DAS}^2$, the number of adaptivity iterations is set to $N_{\rm adaptive}=5$. To demonstrate the effectiveness of the proposed method, we adopt dolfin-adjoint \cite{mitusch2019dolfin} to solve the optimal control problem with some fixed parameters. The dolfin-adjoint solutions, which are regarded as the ground truth, are evaluated on a $200 \times 100$ grid for the physical domain and for $\mb{\xi}$ located on an $11 \times 11$ grid for the parametric domain.

\begin{figure}[!htb]
    \centering
    \subfloat[][The error w.r.t. sample size $|\mathsf{S}_\Omega|$.]{\includegraphics[width=.45\textwidth]{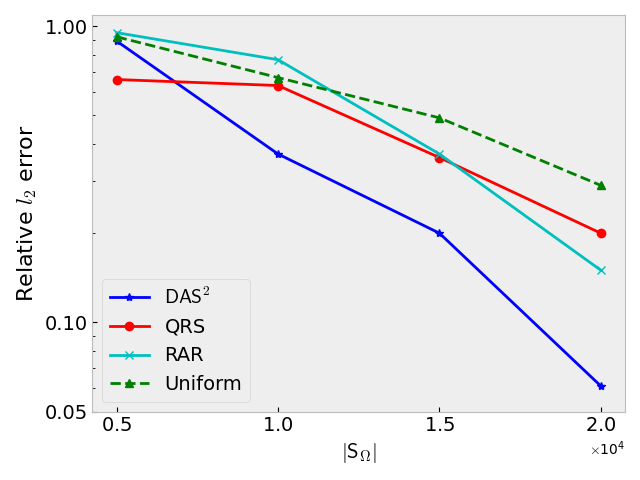}}
    \subfloat[][The error evolution with $|\mathsf{S}_\Omega|=2\times 10^4$.]{\includegraphics[width=.45\textwidth]{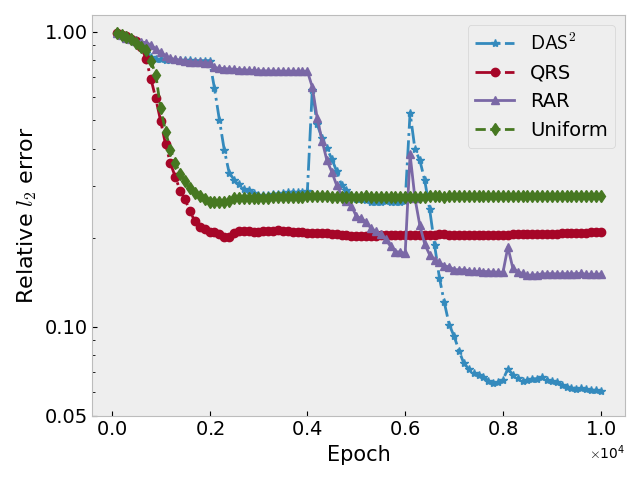}}
    \caption{Approximation errors for the parametric optimal control problem.}
    \label{fig:ex4_comparison}
\end{figure}

\begin{figure}[!htb]
    \centering
    \includegraphics[width=0.55\textwidth]{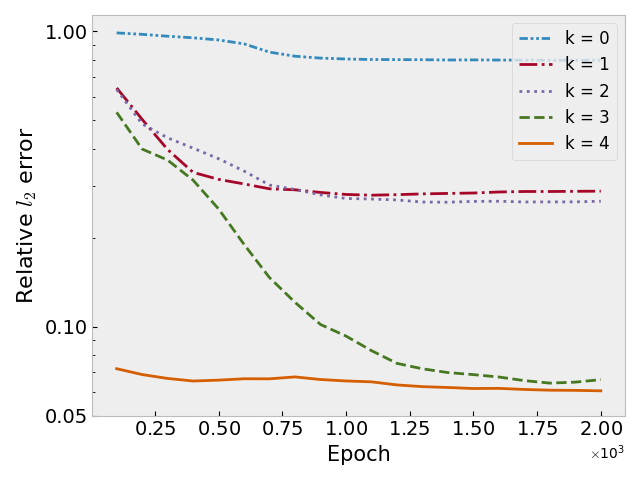}
    \caption{The errors for the parametric optimal control problem of $\mathrm{DAS}^2$ at differen adaptivity iteration steps. $|\mathsf{S}_\Omega|=2\times 10^4.$}
    \label{fig:ex4_das_stage}
\end{figure}

In Figure~\ref{fig:ex4_comparison}, we plot the relative $l_2$ errors given by different sampling strategies with respect to the sample size in the left plot and with respect to the number of epochs in the right plot. For each sample size, we take
three runs with different initialization and compute the mean relative error of the three runs as the final error. For $\mathrm{DAS}^2$, the size of the initial training set $|\mathsf{S}_{\Omega,0}|=n_r$ is set to $1\times10^3, 2\times10^3, 3\times10^3, 4\times10^3$ for $|\mathsf{S}_\Omega|=0.5\times10^4, 1\times10^4, 1.5\times10^4, 2\times10^4$ respectively. For the uniform sampling strategy and the qusi-random sampling (QRS)  strategy, the number of epochs is set to be the same as the total number of epochs of $\mathrm{DAS}^2$, and the number of points in $\mathsf{S}_{\Omega}$ is also set to be the same as $\mathrm{DAS}^2$. For the heuristic method RAR, the numbers of collocation
points in $\mathsf{S}_{\Omega,k}^g$ ($k = 1,2,3,4$) are set to $n_r = 5 \times 10^2,  1\times 10^3, 1.5 \times 10^3, 2 \times 10^3$ for $\vert \mathsf{S}_{\Omega} \vert = 5 \times 10^3, 1 \times 10^4, 1.5 \times 10^4, 2 \times 10^4$ respectively. It is clear that for this test problem $\mathrm{DAS}^2$ has a better performance than the other three (uniform, RAR and QRS) sampling strategies. From the left plot of Figure~\ref{fig:ex4_comparison}, it is clear that, as the number of samples increases, the relative error of $\mathrm{DAS}^2$ decreases faster than those of the uniform sampling strategy, QRS and RAR. Figure~\ref{fig:ex4_comparison}(b) shows the error evolution of different sampling strategies and it is seen that $\mathrm{DAS}^2$ eventually yields a smaller error than the other three sampling methods for the same sample size. Figure~\ref{fig:ex4_das_stage} shows the error evolution of $\mathrm{DAS}^2$ at each adaptivity iteration step $k$. It is seen that, as $k$ increases, the relative error decreases quickly, implying that $\mathrm{DAS}^2$ is effective.

\begin{figure}[!htb]
    \centering
    \subfloat[][ $\mathsf{S}_{\Omega,1}^g$ in RAR]{\includegraphics[width=.33\textwidth]{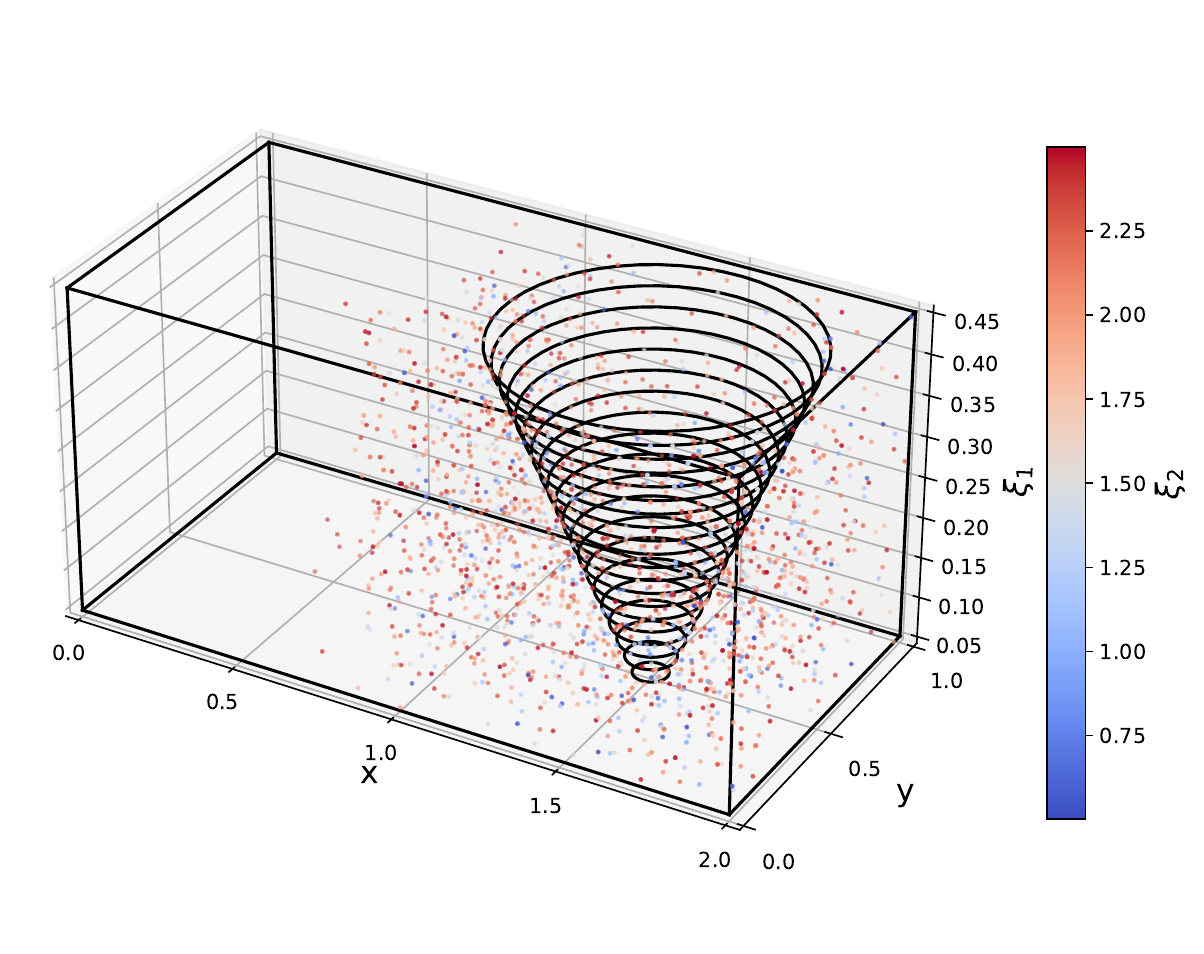}}
    \subfloat[][$\mathsf{S}_{\Omega,2}^g$ in RAR]{\includegraphics[width=.33\textwidth]{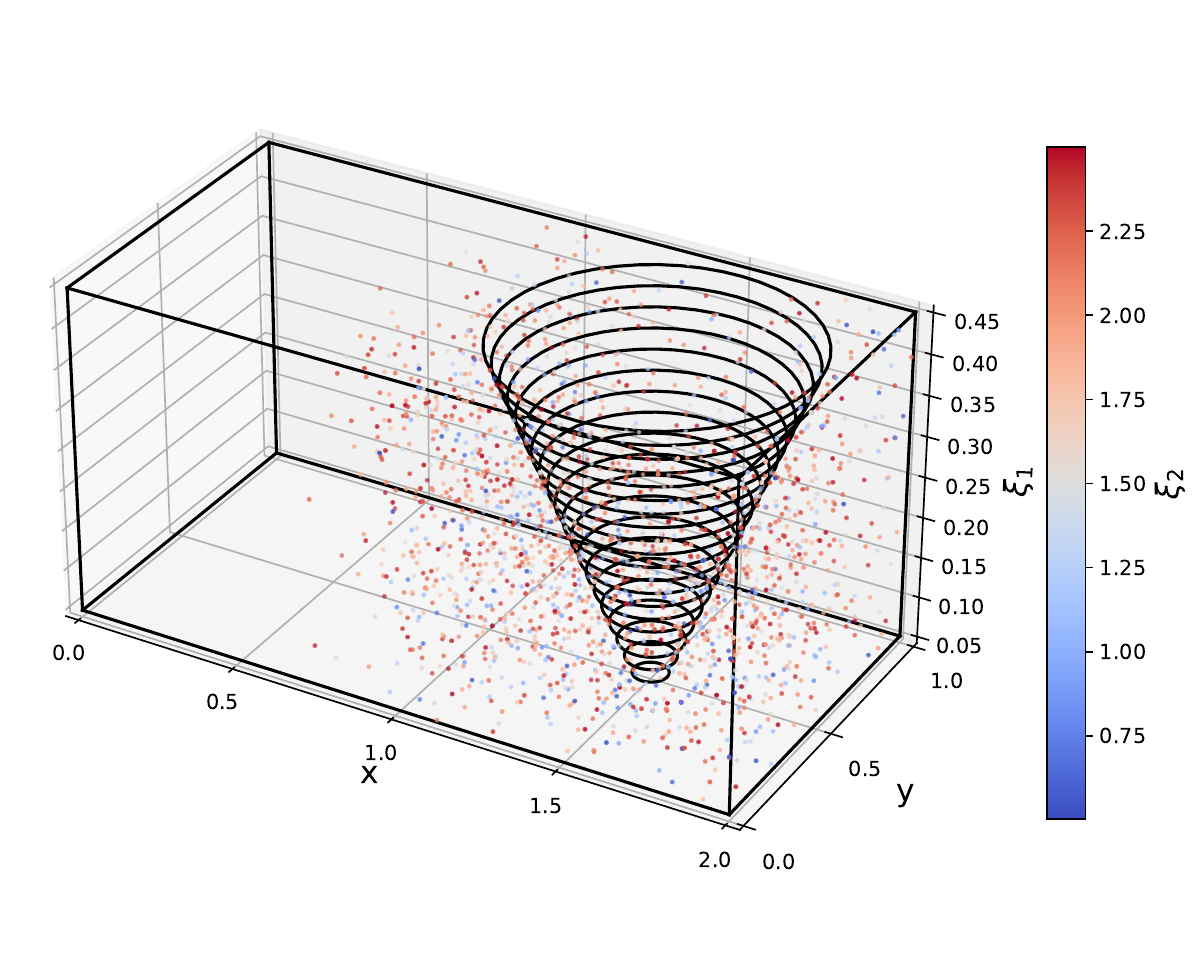}}
    \subfloat[][ $\mathsf{S}_{\Omega,4}^g$ in RAR]{\includegraphics[width=.33\textwidth]{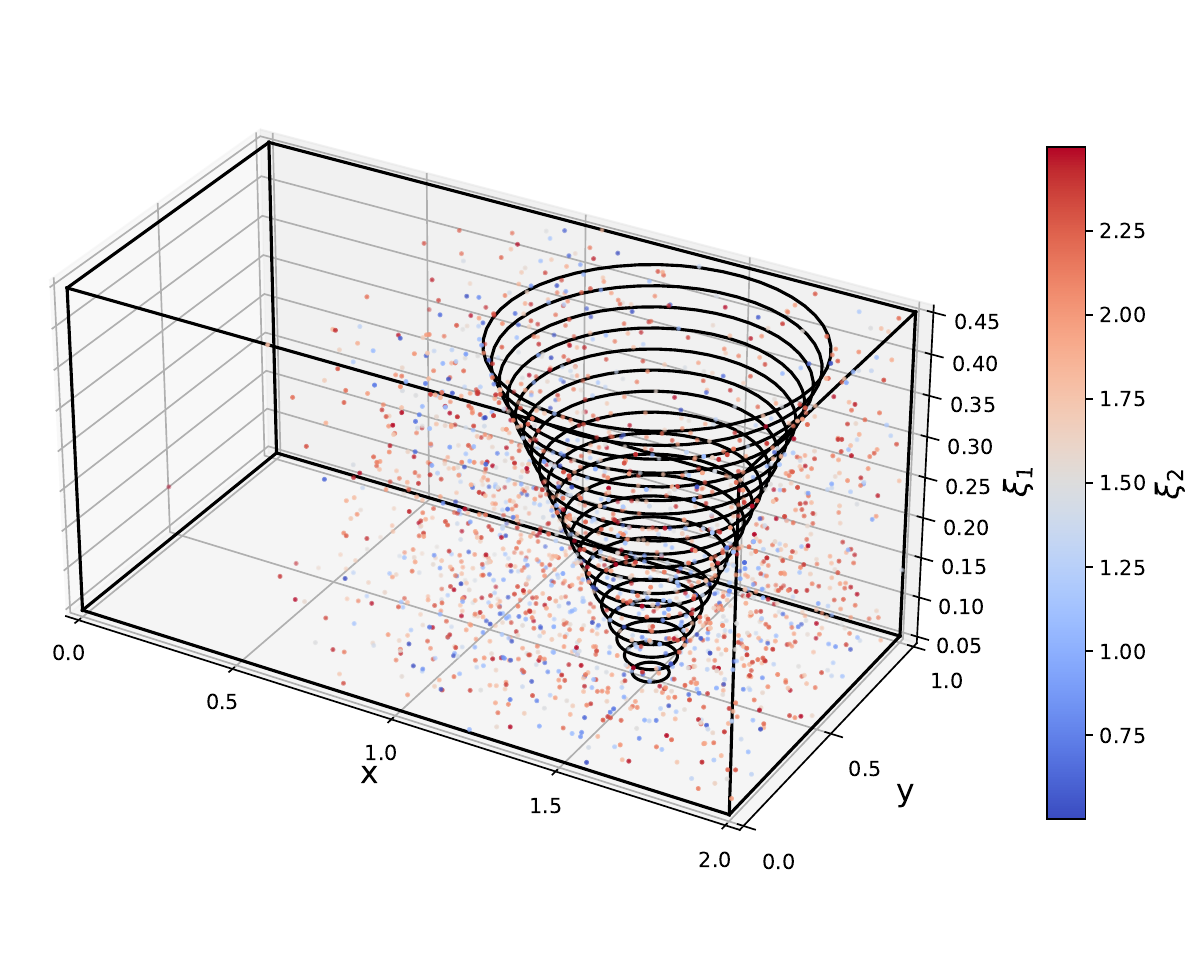}}\\
    \subfloat[][ $\mathsf{S}_{\Omega,1}^g$ in $\mathrm{DAS}^2$]{\includegraphics[width=.33\textwidth]{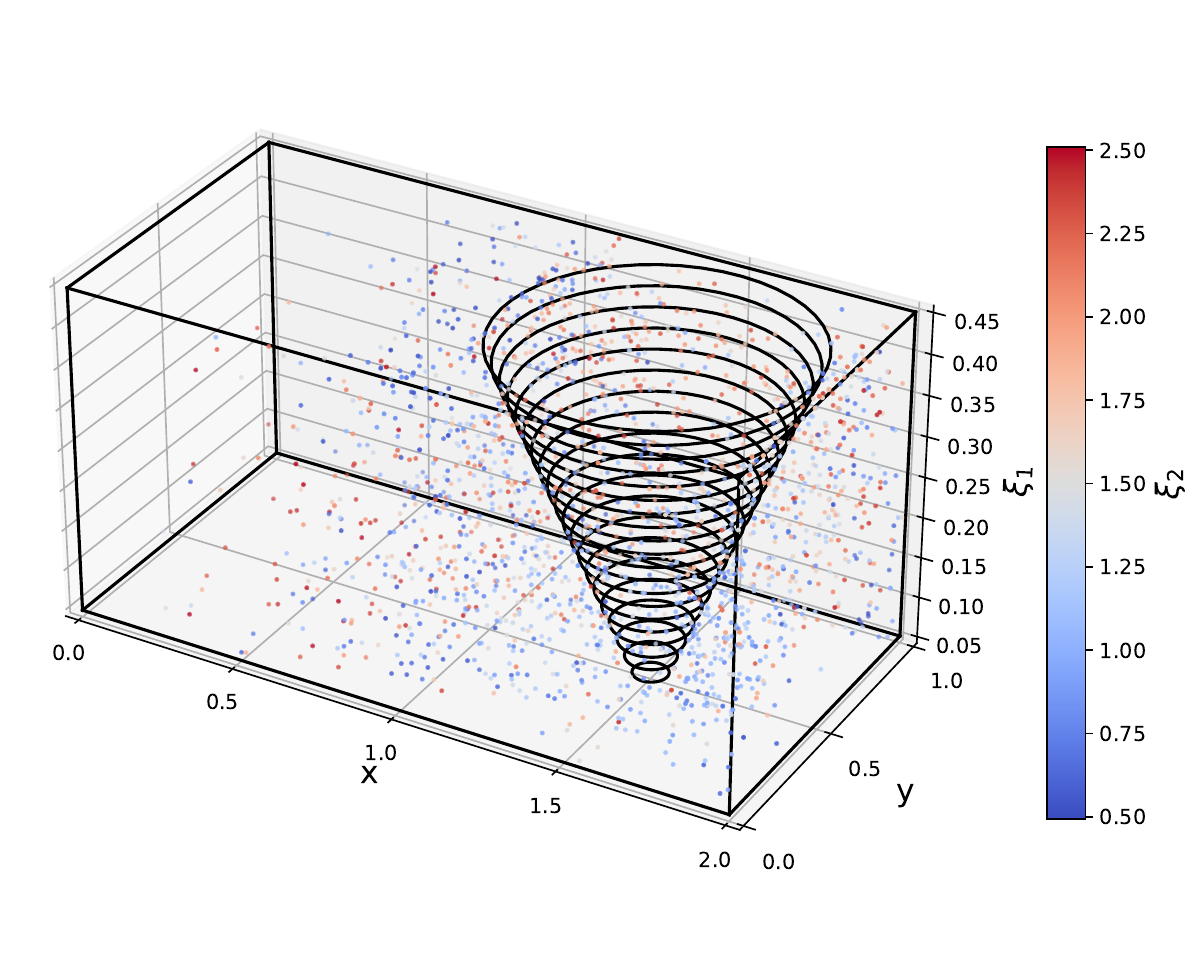}}
    \subfloat[][ $\mathsf{S}_{\Omega,2}^g$ in $\mathrm{DAS}^2$]{\includegraphics[width=.33\textwidth]{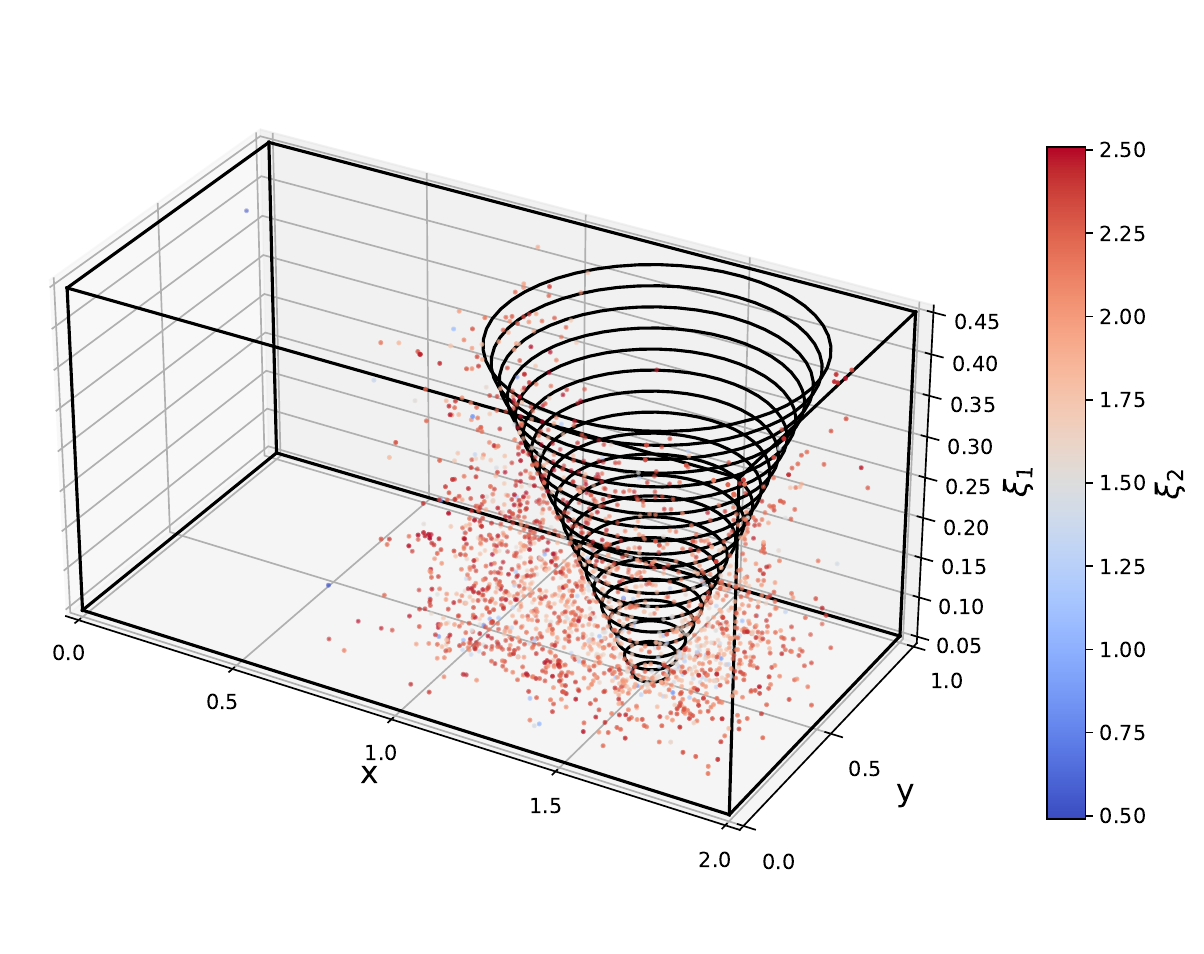}}
    \subfloat[][ $\mathsf{S}_{\Omega,4}^g$ in $\mathrm{DAS}^2$]{\includegraphics[width=.33\textwidth]{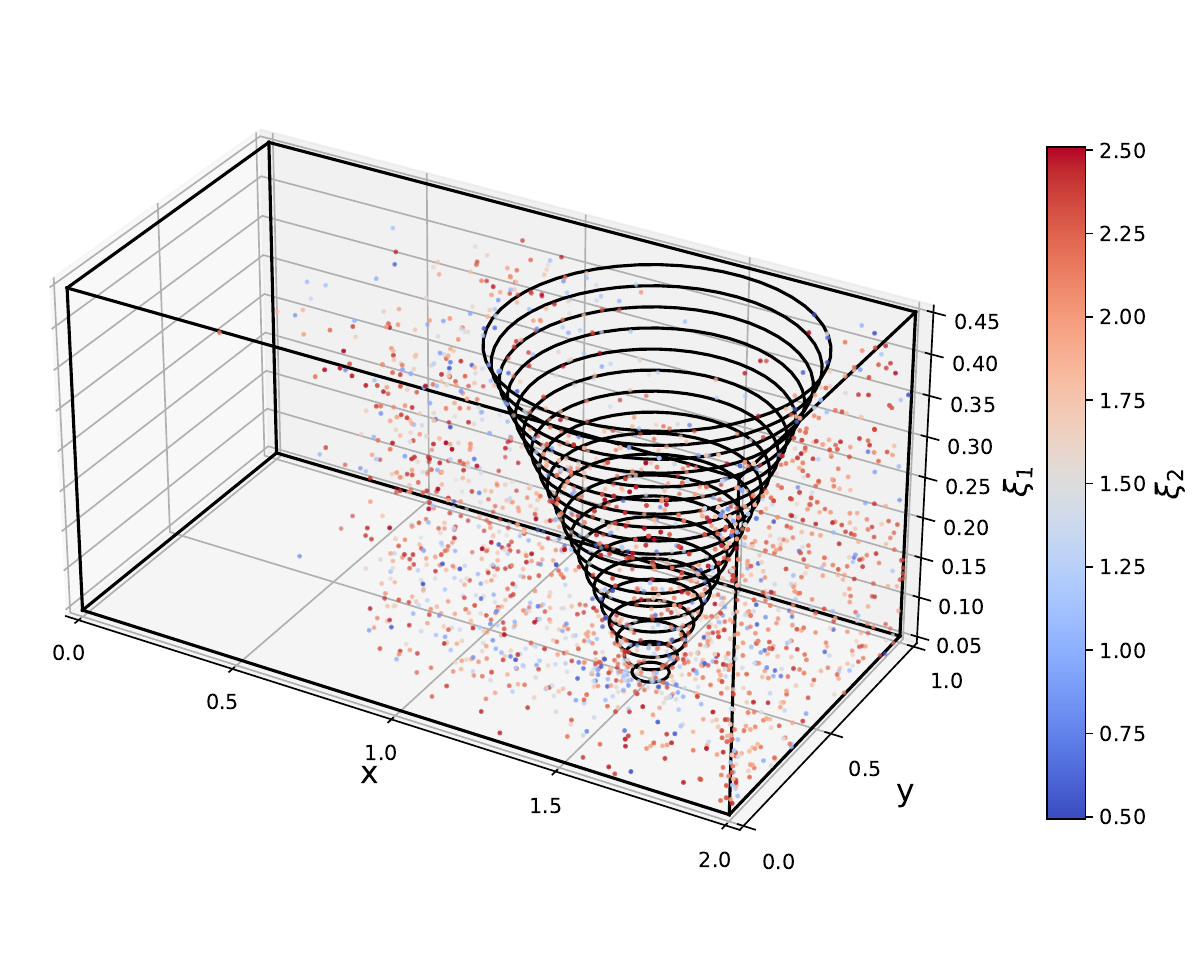}}
    \caption{The evolution of $\mathsf{S}_{\Omega,k}^g$ in RAR and $\mathrm{DAS}^2$ for the parametric optimal control problem (no points inside the frustum), $|\mathsf{S}_\Omega|=2 \times 10^4$.}
    \label{fig:ex4_sampling}
\end{figure}

\begin{figure}[!htb]
    \centering
    \subfloat[][Solutions $u(\mb{x},\mb{\xi})$ with several realizations of $\mb{\xi}$, the parametric optimal control problem.]{\includegraphics[width=.95\textwidth]{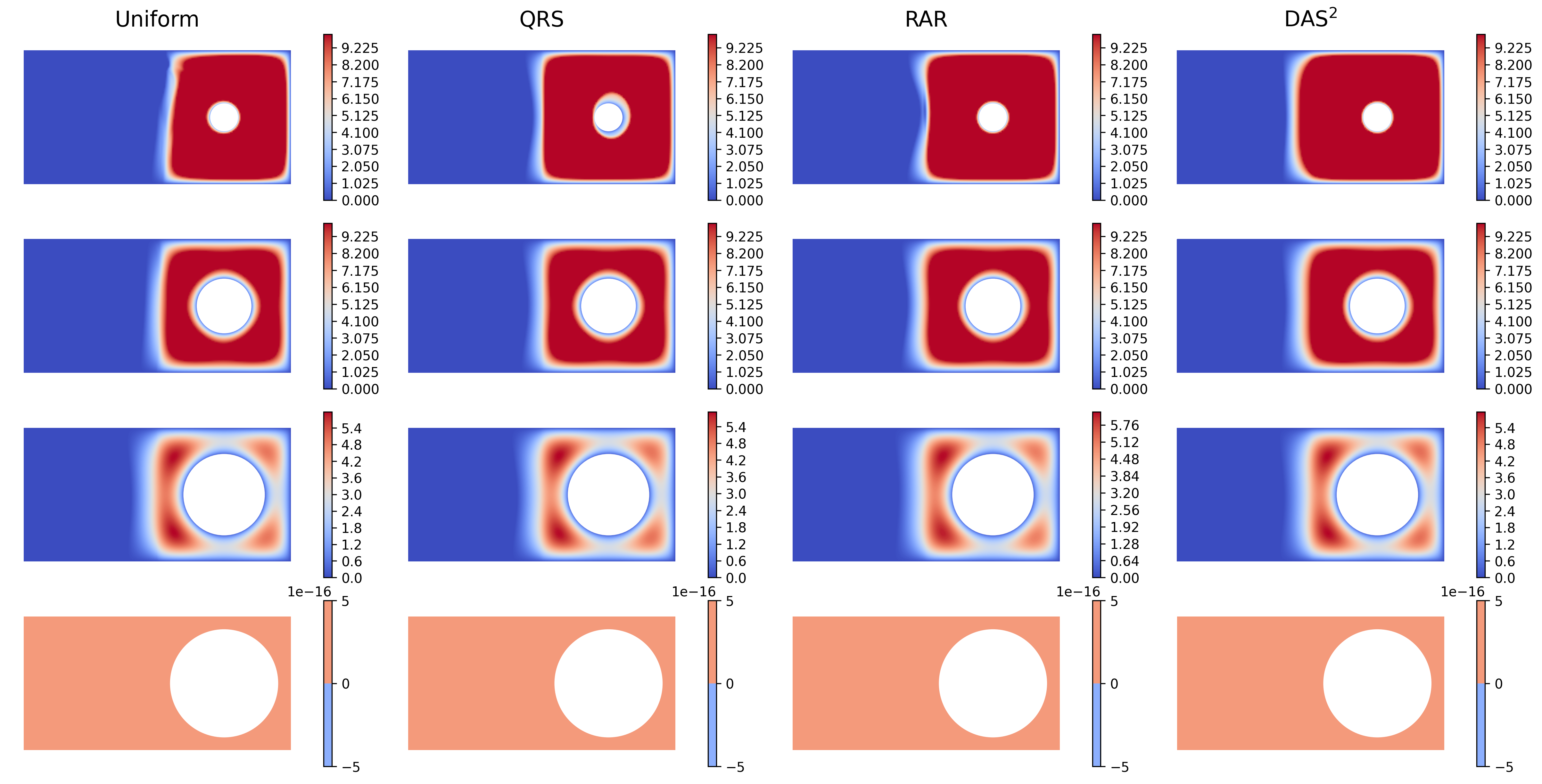}}\\
    \subfloat[][The corresponding absolute point-wise errors for the above solutions, the parametric optimal control problem.]{\includegraphics[width=.95\textwidth]{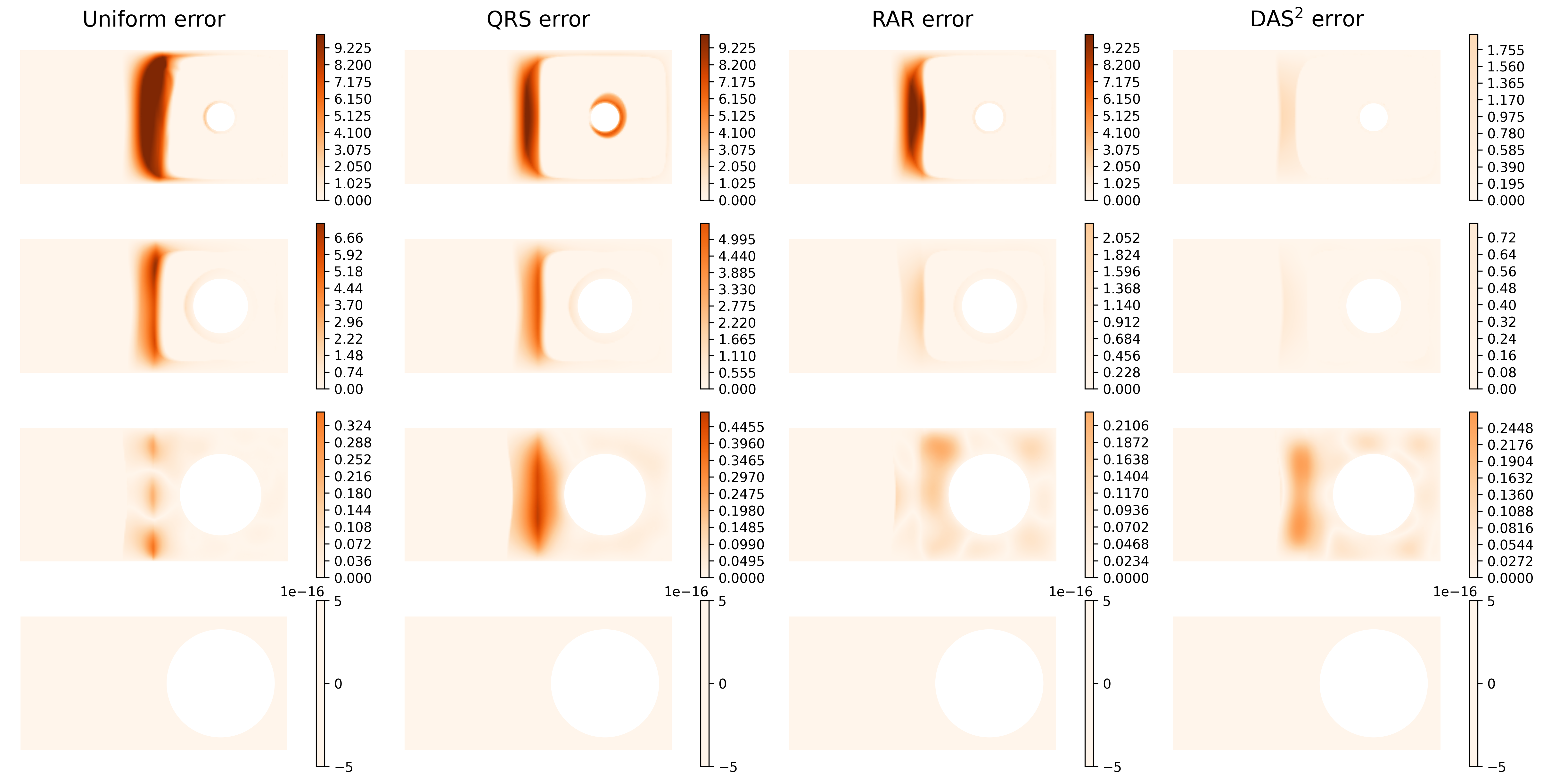}}
    \caption{The results of parametric optimal control problem: the solutions and absolute point wise errors of the uniform sampling method, qusi-random sampling (QRS), RAR and $\mathrm{DAS}^2$ for different realizations. In each subplot, the first line corresponds to $\mb{\xi}=(0.10,2.5)$, the second line corresponds to $\mb{\xi}=(0.20,2.0)$, the third line corresponds to $\mb{\xi}=(0.30,1.5)$, and the fourth line corresponds to $\mb{\xi}=(0.40,0.5)$.}
    \label{fig:ex4_presentation}
\end{figure}

\begin{table}[!htb]
    \caption{The parametric optimal control problem: comparison of different sampling strategies for neural network surrogate models. The relative error of $\mathrm{DAS}^2$ decays faster than other sampling strategies. Once the training of surrogate models is finished, the optimal solution for any parameter can be computed efficiently, which is much faster than the dolfin-adjoint solver (the dolfin-adjoint solver needs 18804 seconds while the neural network surrogate model based on $\mathrm{DAS}^2$ only needs 0.1 seconds).}
    \centering	
    \begin{small}
        \begin{tabular}{cccccccccc}  
            \toprule 
            \multicolumn{2}{c}{\multirow{2}*{\diagbox{sampling strategy}{$\vert \mathsf{S}_{\Omega} \vert$}}} 
            & \multicolumn{1}{c}{$0.5\times 10^4$}& & \multicolumn{1}{c}{$1\times 10^4$} & & \multicolumn{1}{c}{$1.5\times 10^4$}  & & \multicolumn{1}{c}{$2\times 10^4$} \\ 
            \\
            \hline
            \multicolumn{2}{l}{Uniform (0.1s)} &   0.92 & &   0.67 & &   0.49 & & 0.29\\
            \multicolumn{2}{l}{QRS (0.1s)} &    0.66 & &   0.63 & &   0.36& & 0.20 \\
            \multicolumn{2}{l}{RAR (0.1s)} &    0.95& &  0.77 & &   0.37 & & 0.15 \\
            \multicolumn{2}{l}{$\mathrm{DAS}^2$ (0.1s)}    & 0.89 & & 0.37 &  & 0.20 & & 0.06 \\
            \bottomrule
        \end{tabular}
    \end{small}
    \label{table_10d_time_error}
\end{table}

Figure~\ref{fig:ex4_sampling} shows the evolution of the training set ($|\mathsf{S}_\Omega|=2 \times 10^4$) of $\mathrm{DAS}^2$ with respect to adaptivity iterations $k=1,2,4$ (we use $2000$ points in $\mathsf{S}_{\Omega,k}^g$ for visualization), where the initial training set $\mathsf{S}_{\Omega,0}$ consists of uniform samples on $\Omega$. Note that $\xi_1$ denotes the radius of the circle and $\xi_2$ is the desired state in $\Omega_2$. We use different colors to identify $\xi_2$ in Figure~\ref{fig:ex4_sampling}. It can be seen that $\mathrm{DAS}^2$ can effectively capture the information of singularity since the data points generated by $\mathrm{DAS}^2$ are concentrated on the area where large residuals are located (see $\mathsf{S}_{\Omega,1}^g$ and $\mathsf{S}_{\Omega,2}^g$ in $\mathrm{DAS}^2$), while RAR is not able to capture the variation in residual well enough. 
Finally, nearly uniform samples are generated to augment the training set in $\mathrm{DAS}^2$ because one can obtain a flat residual profile after four adaptivity iterations. Figure~\ref{fig:ex4_presentation}(a) shows the optimal control solution $u(\mb{x},\mb{\xi})$ obtained using different sampling methods. We choose several different parameters $\mb{\xi}$ for visualization. For validation, the absolute errors between different sampling strategies and the dolfin-adjoint solver are plotted in Figure~\ref{fig:ex4_presentation}(b). It can be seen that the $\mathrm{DAS}^2$ has a better performance than the other three sampling strategies. Table~\ref{table_10d_time_error} shows the inference time and the relative error for the uniform sampling strategy, QRS, RAR and $\mathrm{DAS}^2$. It is seen that $\mathrm{DAS}^2$ performs much better than the other three sampling strategies especially when the sample size is relatively large. 

	
\subsection{Surrogate modeling for parametric lid-driven cavity flow problems with $Re \in [100, 1000]$}\label{sec_ldc_flow}
Finally, we consider the lid-driven cavity flow problem governed by the following  steady-state incompressible Navier-Stokes equations 
\begin{equation}
    \left\{
    \begin{aligned}
        \mb{u}(\mb{x},\xi)\cdot \nabla \mb{u}(\mb{x},\xi) + \nabla p(\mb{x},\xi) = \frac{1}{Re(\xi)} \Delta \mb{u}(\mb{x},\xi) \quad & \text{in}\;\Omega, \\
        \nabla \cdot \mb{u}(\mb{x},\xi)=0 \quad & \text{in}\; \Omega, \\
        \mb{u}(\mb{x},\xi) = \mb{g}(\mb{x},\xi) \quad & \text{on}\;\partial \Omega,
    \end{aligned}
    \right.
    \label{eq:lid-driven cavity}
\end{equation}
where $\mb{u}(\mb{x},\xi) = [u(\mb{x},\xi), v(\mb{x},\xi)]^{\mathsf{T}}$ and $p(\mb{x},\xi)$ are the flow velocity field and the scalar pressure respectively. Here, we consider a parametric problem in terms of the Reynolds number, where we assume that  
$Re(\xi) = \xi \in \Omega_p = [100,1000]$. The physical domain is $\Omega_s = [0,1] \times [0,1]$. The velocity profile $\mb{u} = [1,0]^{\mathsf{T}}$ is imposed on the top boundary ($y = 1$ where $\mb{x} = [x, y]^{\mathsf{T}}$), and $\mb{u} = [0,0]^{\mathsf{T}}$ is imposed on all other boundaries, i.e., for $[\mb{x}, \xi] \in \partial \Omega $
\begin{equation*}
    \mb{g}(\mb{x},\xi) = \begin{cases}
        [1,0]^{\mathsf{T}}, y = 1;\\
        [0,0]^{\mathsf{T}}, \ \text{otherwise}.
    \end{cases}
\end{equation*}
The lid-driven cavity problem is a benchmark in computational fluid dynamics. However, even for a fixed relatively low Reynolds number, the existing neural-network-based methods are not able to achieve a comparable accuracy with the baseline obtained by classical numerical methods \cite{ghia1982high}. In this study, we use the proposed $\mathrm{DAS}^2$ method to obtain accurate all-at-once solutions of the parametric lid-driven cavity flow problem with Reynolds numbers from $100$ to $1000$.

We first evaluate the performance of $\mathrm{DAS}^2$ with the non-parametric lid-driven cavity flow problem, where we consider $Re=100,400$. In such a scenario, $\mathrm{DAS}^2$ reduces to $\mathrm{DAS}$. After that, we use $\mathrm{DAS}^2$ to solve the parametric lid-driven cavity flow problem to obtain all-at-once solutions, where we consider $Re\in[100,1000]$. To measure the quality of $\mathrm{DAS}^2$, we compare $\mathrm{DAS}^2$ with the classical numerical methods presented in the literature \cite{ghia1982high}, the FEniCS solver \cite{alnaes2015fenics,logg2012automated}, and the WAM-AW method proposed in a recent literature \cite{hou2023enhancing}. 
Unlike DAS, WAM-AW is an adaptive collocation point movement approach based on interacting particle methods for solving low-regularity PDEs, and we also use this method as a baseline of neural-network-based methods.

Since there are three quantities ($u,v$ and $p$) to be determined in equation \eqref{eq:lid-driven cavity}, we construct a neural network $\mb{u}_{\mb{\theta}}$ with three outputs to represent $u,v$ and $p$. For the deterministic problem (fixed Reynolds numbers), we choose a five-layer fully connected neural network $\mb{u}_{\mb{\theta}}(\mb{x})$, where each hidden layer has 20 neurons.  
For KRnet, we set $K = 2$ and the configuration for the affine coupling layers remains the same as the previous experiments. 
The number of epochs for training both $\mb{u}_{\mb{\theta}}(\mb{x})$ and $p_{\mathsf{KRnet}}(\mb{x};\mb{\theta}_f)$ is set to $N_e = 3000$. The optimizer for training $\mb{u}_{\mb{\theta}}(\mb{x})$ is BFGS, the optimizer for training $p_{\mathsf{KRnet}}(\mb{x};\mb{\theta}_f)$ is ADAM with learning rate $0.0001$ and the batch size is set to $m = 100$ (for $Re = 100$) or $m=500$ (for $Re = 400$). For $Re = 100$, the number of adaptivity iterations is set to $N_{\rm adaptive}=5$ with $n_r=200$, resulting in $|\mathsf{S}_{\Omega}|=1000$. For the boundary term, $400$ boundary points are uniformly sampled on $\partial \Omega_s$ with $100$ points for each edge. For the WAM-AW method proposed in the literature \cite{hou2023enhancing}, we exactly keep the setup of their work and run their open source code, where the number of collocation points is also set to $|\mathsf{S}_{\Omega}|=1000$ and $100$ data points for each edge of boundary (these settings are the same as in \cite{hou2023enhancing}). For $Re = 400$, the number of adaptivity iterations is set to $N_{\rm adaptive}=10$ with $n_r=500$, resulting in $|\mathsf{S}_{\Omega}|=5000$. We uniformly sample $256$ points for each edge on the boundary. For the WAM-AW method, we also set $|\mathsf{S}_{\Omega}|=5000$ and 256 boundary points are sampled on each edge of the boundary. 
For all cases, we discretize in space using the $Q_3$-$Q_2$ finite element method implemented in FEniCS with a uniform $129\times 129$ grid to obtain a reference solution. 

\begin{figure}[!htb]
    \centering
    \subfloat[][$Re=100.$]{\includegraphics[width=.45\textwidth]{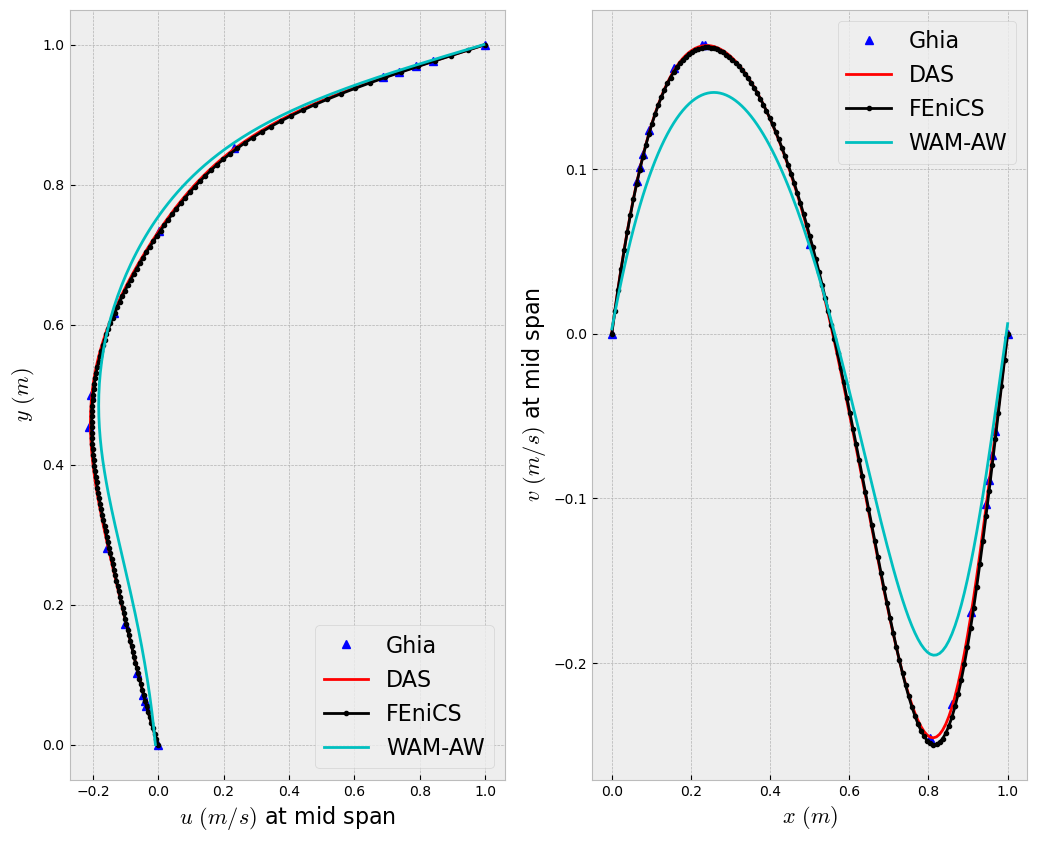}}
    \subfloat[][$Re=400$.]{\includegraphics[width=.45\textwidth]{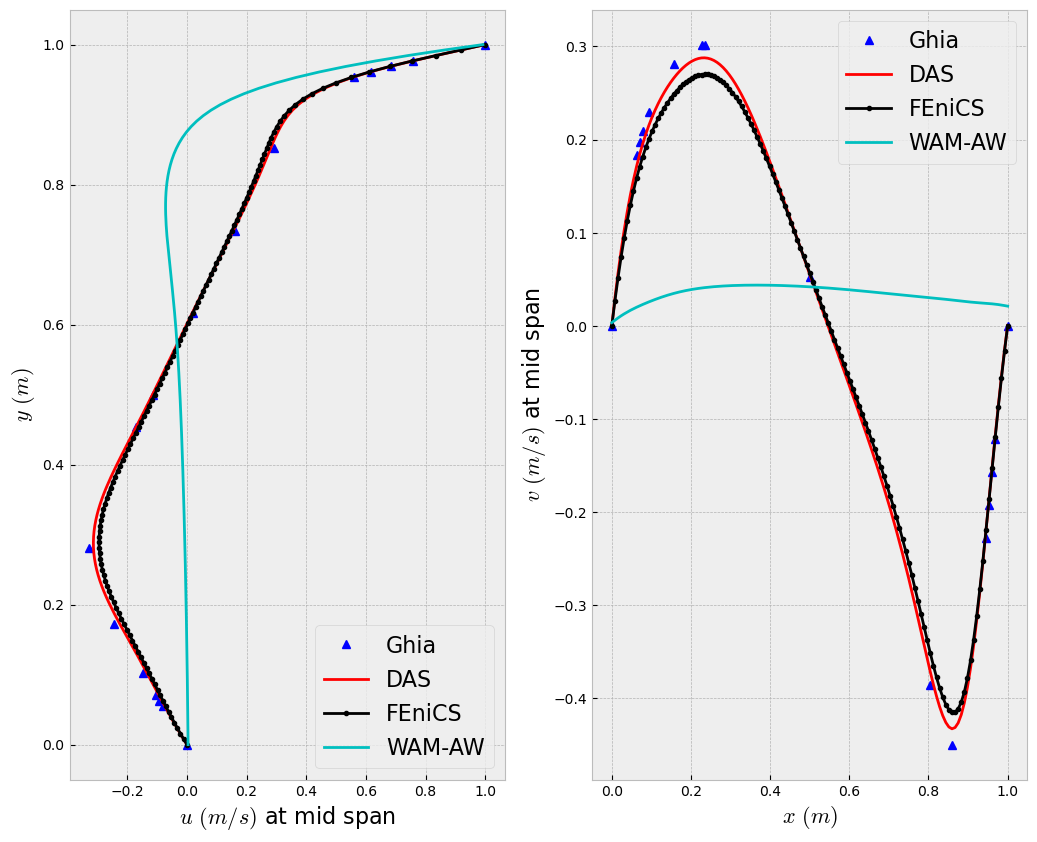}}
    \caption{The velocity components at the location of mid-span lines for the deterministic lid-driven cavity flow problems, $Re = 100, 400$.}
    \label{fig:ldc_fixed_Re_comparison}
\end{figure}

\begin{figure}[!htb]
    \centering  \includegraphics[width=.8\textwidth]{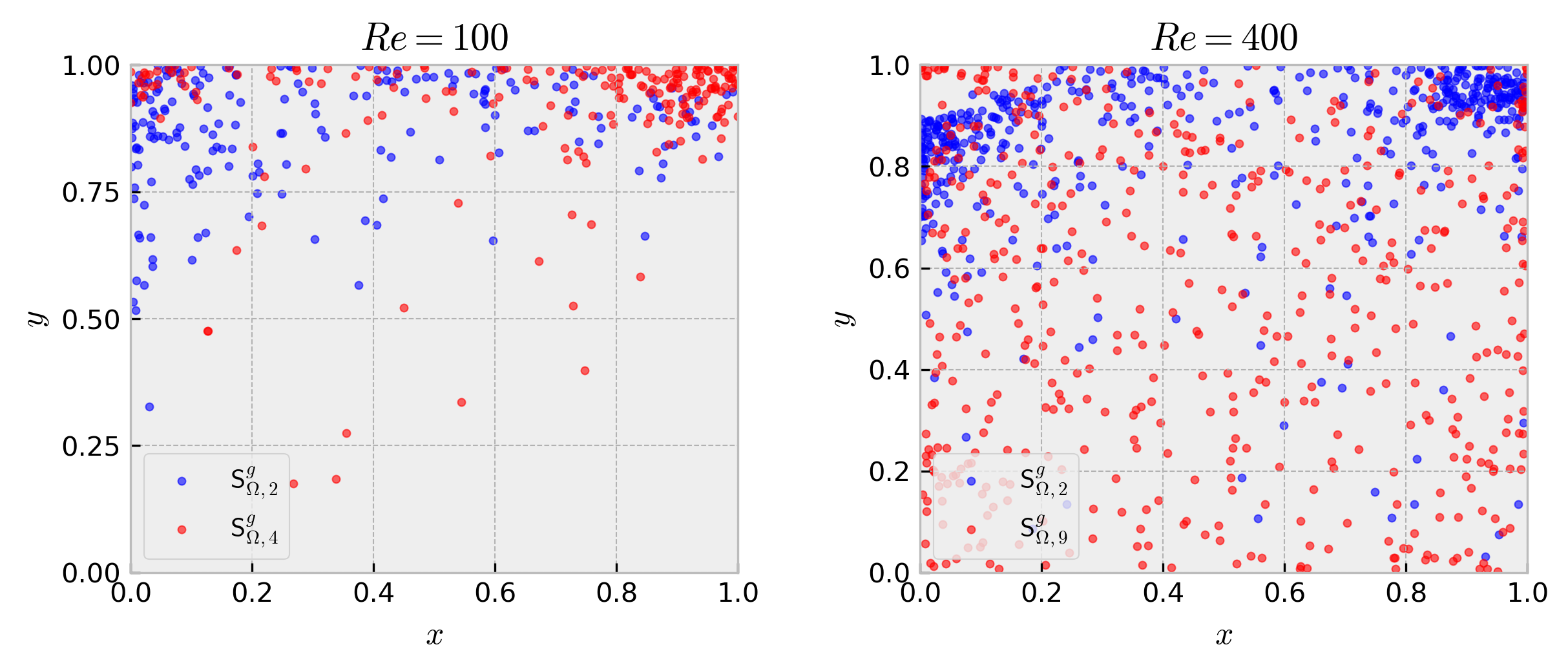}
    \caption{The random samples in $\mathsf{S}_{\Omega, k}^g$ for the deterministic lid-driven cavity flow problems. Left: $\mathsf{S}_{\Omega, 2}^g$ (blue) and $\mathsf{S}_{\Omega, 4}^g$ (red) for $Re = 100$; Right: $\mathsf{S}_{\Omega, 2}^g$ (blue) and $\mathsf{S}_{\Omega, 9}^g$ (red) for $Re = 400$.}
    \label{fig:ex2_fixed_Re_samples}
\end{figure}

\begin{figure}[!htb]
    \centering
    \subfloat[][$Re=100.$]{\includegraphics[width=0.85\textwidth]{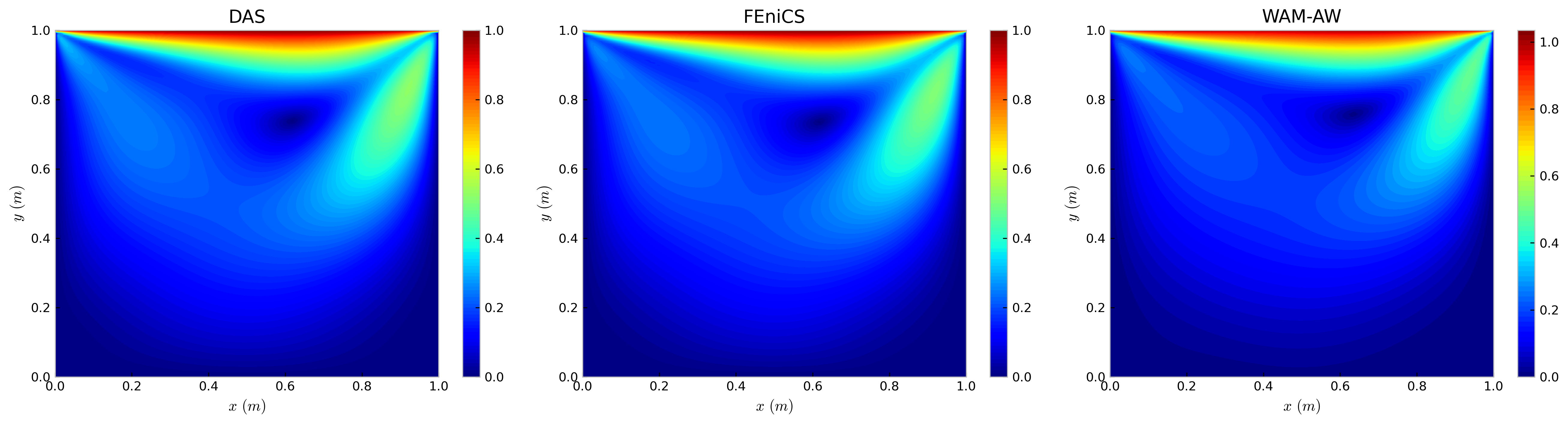}}\\
    \subfloat[][$Re=400.$]{\includegraphics[width=0.85\textwidth]{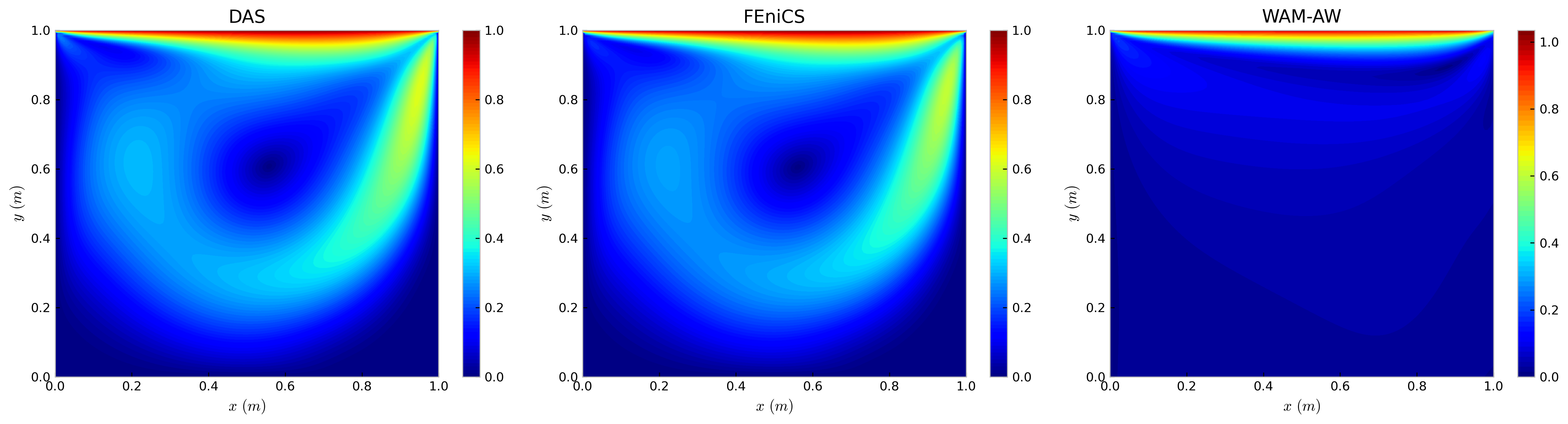}}
    \caption{The visualization of $|\mb{u}| = \sqrt{u^2 + v^2}$ for the deterministic lid-driven cavity flow problems.}
    \label{fig:ldc_fixed_Re_sol}
\end{figure}
	
Figure~\ref{fig:ldc_fixed_Re_comparison} shows the velocity at the location of the mid-span line, which is usually used to assess the accuracy of solutions. Specifically, for $x = 0.5$, we plot the velocity component $u$ with respect to $y$ and for $y=0.5$ we plot the velocity component $v$ with respect to $x$. In Figure~\ref{fig:ldc_fixed_Re_comparison}, we compare the results of DAS, FEniCS, WAM-AW with the benchmark results given in Ghia et. al \cite{ghia1982high}. It is seen that the results given by DAS are consistent with those given by Ghia and FEniCS, while the results of WAM-AW do not agree with the reference results, especially for $Re = 400$. 
To further illustrate the effectiveness of DAS, we plot the evolution of random samples during training in Figure~\ref{fig:ex2_fixed_Re_samples}, where the left plot shows $\mathsf{S}_{\Omega,2}^g$ and $\mathsf{S}_{\Omega,4}^g$ for $Re = 100$ and the right plot shows $\mathsf{S}_{\Omega,2}^g$ and $\mathsf{S}_{\Omega,9}^g$ for $Re = 400$. It can be seen that DAS yields samples that are consistent with both the problem properties and the approximation, where the initial training set consists of random samples generated by Latin hypercube sampling.  
For example, at $k = 2$, most of the samples in $\mathsf{S}_{\Omega,2}^g$ are located in the upper corners, where the velocity field changes abruptly and large residuals occur. As $k$ increases, the residual profile becomes more uniform after the localized information is well captured, which implies that random samples can be added more uniformly. 
Figure~\ref{fig:ldc_fixed_Re_sol} shows the image of $\vert \mb{u} \vert = \sqrt{u^2 + v^2}$, where $Re = 100$ and $Re = 400$ are considered. Compared with the reference solution given by FEniCS, DAS provides an accurate prediction of the flow velocity for $Re=100$ while WAM-AW has a little loss of accuracy. For $Re=400$, the results given by DAS are still accurate while the results given by WAM-AW are not physically correct. 

Next, we look at the surrogate modeling of parametric lid-driven cavity flow problems. The architecture of neural networks and the setting of adaptive sampling for surrogate modeling need to be modified since solving such parametric problems is more difficult than deterministic ones. We use one five-layer fully connected neural network $\mb{u}_{\mb{\theta}}(\mb{x},\xi)$ with three outputs to approximate the parametric solutions $u(x,y,\xi)$, $v(x,y,\xi)$, $p(x,y,\xi)$ respectively, where each hidden layer has $32$ neurons.  For adaptive sampling, we use the joint PDF model induced by KRnet in $\mathrm{DAS}^2$. For KRnet, we set $K=3$ and take $L=6$ affine coupling layers. For each affine coupling layer, a two-layer fully connected neural network is employed where each hidden layer has $32$ neurons. The number of epochs for training the surrogate model and KRnet is set to $N_e=5000$. The optimizer for training the surrogate model $\mb{u}_{\mb{\theta}}$ is BFGS, and the optimizer for training KRnet is ADAM with a learning rate $0.0001$.  
The number of adaptivity iterations is set to $N_{\rm adaptive}=10$ with $n_r=1 \times 10^4$, resulting in the total number of collocation points $|\mathsf{S}_{\Omega}|=1 \times 10^5$. For the boundary term, $16384$ boundary points are sampled on each edge of the boundary. The batch size is set to $m=5000$. 

\begin{figure}[!htb]
    \centering
    \includegraphics[width=1.\textwidth]{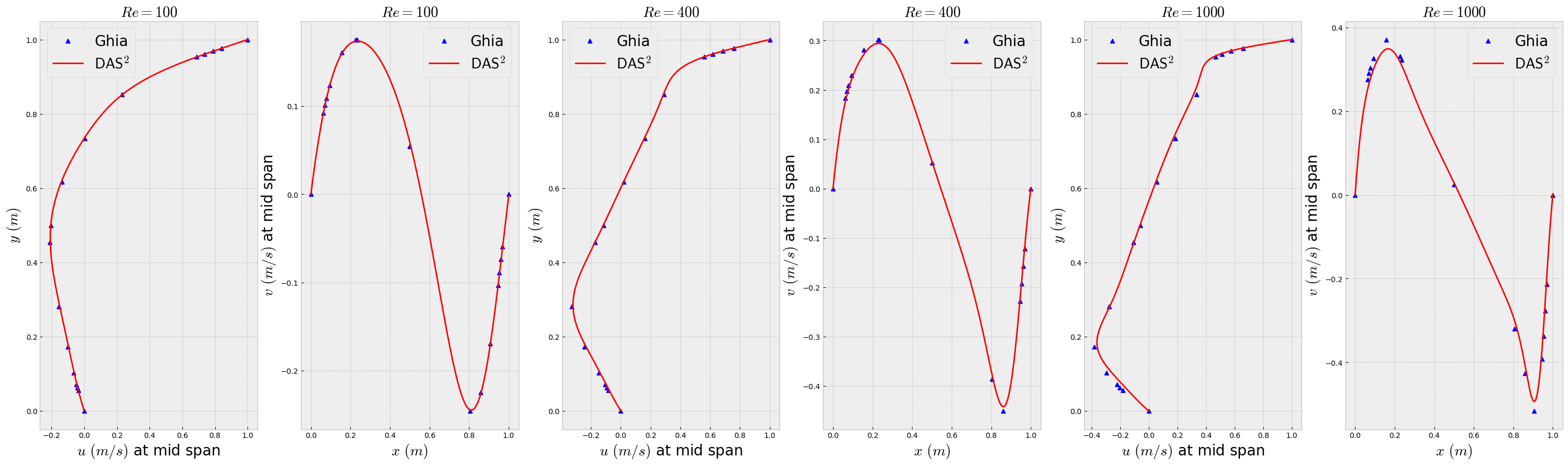}
    \caption{The velocity components at the location of mid-span lines for surrogate modeling of parametric lid-driven cavity flow problems ($Re \in [100,1000]$). The results for $Re = 100, 400, 1000$ are chosen for visualization.}
    \label{fig:ex2_ghia_all_at_once}
\end{figure}

\begin{figure}[!htb]
    \centering
   {\includegraphics[width=0.75\textwidth]{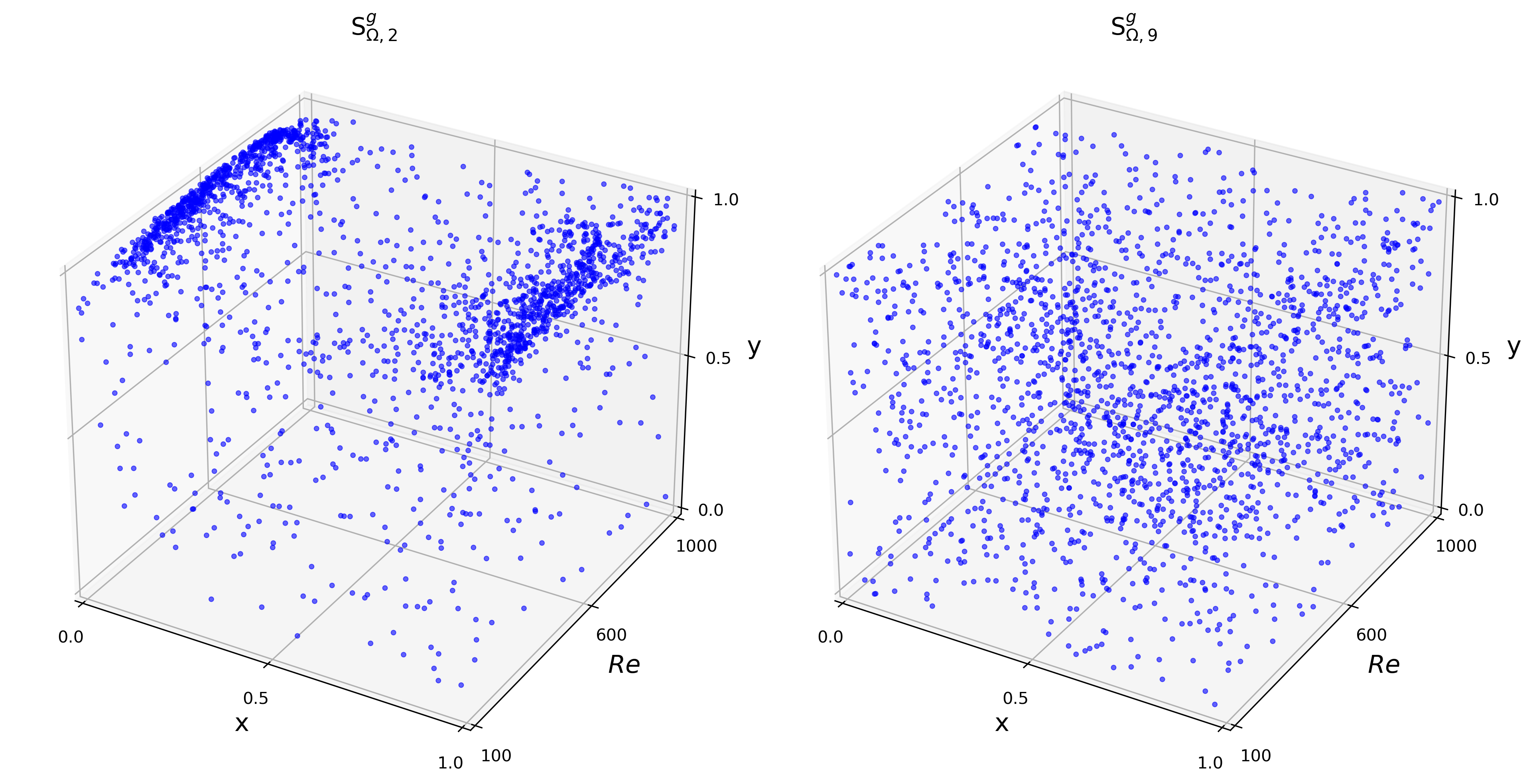}}\\
    \caption{The random samples in $\mathsf{S}_{\Omega, k}^g$ (2000 samples are displayed) for surrogate modeling of parametric lid-driven cavity flow problems, $Re \in [100,1000]$.}
    \label{fig:ldc_das2_sample}
\end{figure}

\begin{figure}[!htb]
    \centering
    \includegraphics[width=1.\textwidth]{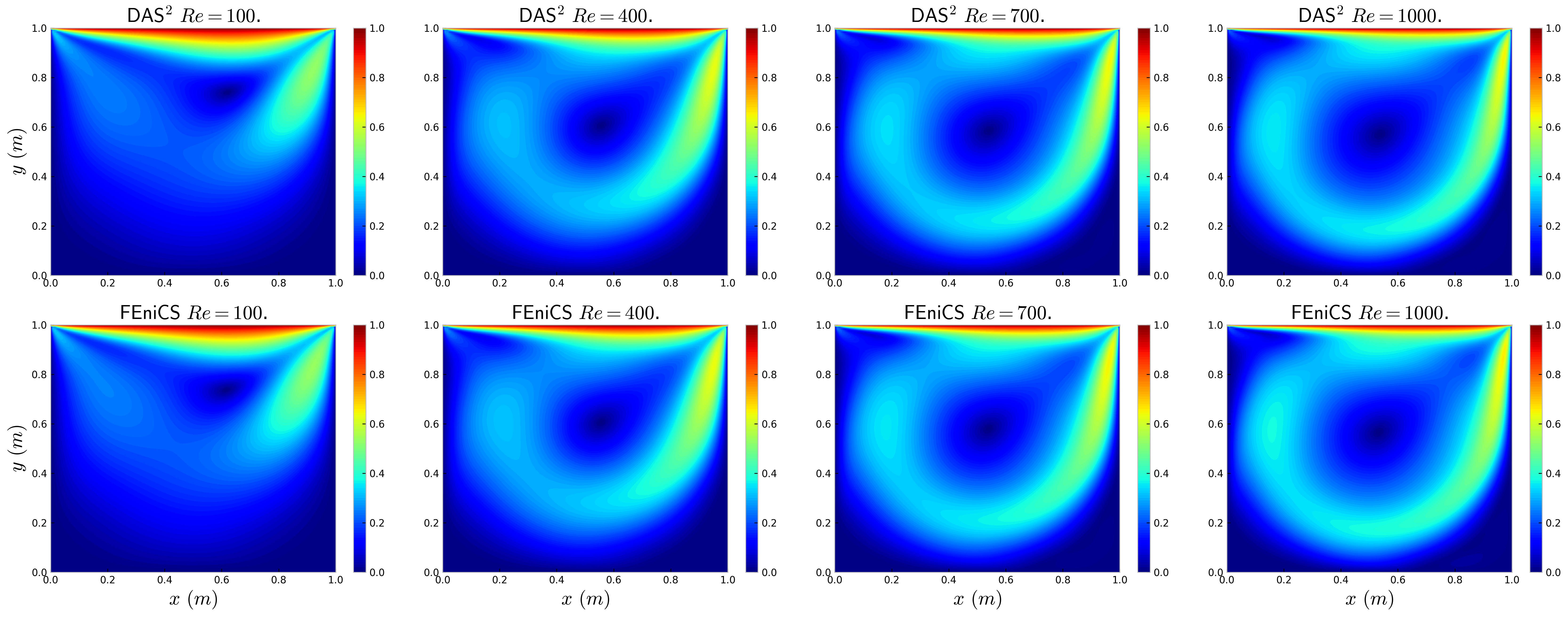}
    \caption{The visualization of $|\mb{u}| = \sqrt{u^2 + v^2}$ for surrogate modeling of parametric lid-driven cavity flow problems, $Re \in [100,1000]$. The $l_2$ relative errors are $1.5\%, 1.1\%, 3.1\%, 4.8\%$ for $Re = 100, 400, 700, 1000$ respectively.}
    \label{fig:ex2_all_at_once}
\end{figure}

Figure~\ref{fig:ex2_ghia_all_at_once} shows the velocity profile given by the trained surrogate model $\mb{u}_{\mb{\theta}}(\mb{x},\xi)$ at the location of the mid-span line for some selected Reynolds numbers $\xi = 100, 400, 1000$.  
From Figure~\ref{fig:ex2_ghia_all_at_once}, it is clear that the results of $\mathrm{DAS}^2$ are consistent with Ghia's data \cite{ghia1982high}, implying that our $\mathrm{DAS}^2$ approach is able to provide an accurate surrogate model for fast inference. Figure~\ref{fig:ldc_das2_sample} shows the evolution of the training set ($|\mathsf{S}_\Omega|=10^5$) of $\mathrm{DAS}^2$ with respect to adaptivity iterations $k=2,9$, where the initial training set $\mathsf{S}_{\Omega,0}$ consists of random samples generated by Latin hypercube sampling. $\mathsf{S}_{\Omega, 2}^g$ indicates that the residual concentrates on the two upper corners for any $Re \in [100,1000]$. As the adaptivity iteration $k$ increases, the residual profile becomes more flat as shown by the distribution of $\mathsf{S}_{\Omega, 9}^g$, which is expected since more collocation points are added to the two upper corners to reduce the errors over there. Figure~\ref{fig:ex2_all_at_once} shows the image of $\vert \mb{u} \vert = \sqrt{u^2 + v^2}$, where $Re = 100, 400, 700, 1000$ are used for visualization. Here, we again use the $Q_3$-$Q_2$ finite element method implemented in FEniCS with
a uniform $129 \times 129$ grid to obtain the reference solutions for $Re = 100, 400, 700, 1000$. It is seen that 
$\mathrm{DAS}^2$ provides an accurate prediction of the flow velocity even for $Re = 1000$. 
The $l_2$ relative errors, which are evaluated on the $129\times129$ uniform grid, are $1.5\%, 1.1\%, 3.1\%, 4.8\%$ for $Re = 100, 400, 700, 1000$ respectively. It is worth noting that the inference time of $\mathrm{DAS}^2$ is 0.02 seconds, while the computation time of FEniCS is 309.94 seconds to obtain the four solutions for $Re = 100, 400, 700, 1000$.


	\section{Conclusions}
	In this paper, we have developed a deep adaptive sampling approach for surrogate modeling ($\mathrm{DAS}^2$) of parametric differential equations, generalizing the previous work DAS to the parametric setting. It has been shown that $\mathrm{DAS}^2$ not only provides a fast inference for parametric differential equations without labeled data but also yields an accurate prediction for low-regularity problems thanks to the adaptive sampling procedure. Similar to DAS, the framework of $\mathrm{DAS}^2$ also utilizes a deep generative model to generate collocation points that are consistent with the residual-induced distribution. Unlike DAS, $\mathrm{DAS}^2$ handles the low regularity from both spatial and parametric spaces. The joint PDF (marginal PDF) model for both spatial and parametric variables (only the parametric variable), which is induced by the residual, provides effective samples to reduce the statistical errors from the discretization of the loss functional. Based on this, the accuracy of all-at-once solutions of parametric differential equations can be improved significantly. 
	
	We pay particular attention to the following observations: First, due to the physics-informed model with adaptive sampling, we have not used any simulation data for the training process. Second, the procedure of adaptive sampling is independent of the structure of the surrogate model unless the model is defined on a certain set of collocation points. Third, deep generative modeling plays an important role in $\mathrm{DAS}^2$. Deep generative modeling outperforms classical density models or sampling strategies in the sense that it effectively merges density approximation and sample generation for an arbitrary high-dimensional distribution. $\mathrm{DAS}^2$ will find many applications because it provides a general way to improve the training set and any improvement in the model structure can be further refined by adaptive sampling. 
	

	\bigskip
	\textbf{Acknowledgments:}
	K. Tang has been supported by the China Postdoctoral Science Foundation grant 2022M711730. J. Zhai is supported by the start-up fund of ShanghaiTech University (2022F0303-000-11). X. Wan has been supported by NSF grant DMS-1913163. C. Yang has been supported by NSFC grant 12131002 and Huawei Technologies Co., Ltd.

 	\appendix
	\begin{appendices}
		\section{Proof of Theorem \ref{thm_bounded_by_residual}}
		\begin{proof}
			Since $\Theta$ is compact, there exists a $\delta$-net $\bar{\mathsf{N}} = \{\mb{\theta}_1, \ldots, \mb{\theta}_{N_{\delta}} \}$ with the following property: for all $\mb{\theta} \in \Theta$, there exists $1 \leq i \leq N_{\delta}$ such that $\norm{\mb{\theta} - \mb{\theta}_i}{\infty} \leq \delta$ \cite{wright2021high}. For a given $\varepsilon \in (0,1)$, we set $\delta = \varepsilon^2/(4 \mathfrak{L})$. Moreover, the number of parameters of $\bar{\mathsf{N}}$ is at most $(4a \mathfrak{L}/\varepsilon^2)^D$. By Assumption \ref{assump_lip_operator}, for $\mb{\theta}, \mb{\nu} \in \Theta$ with $\norm{\mb{\theta} - \mb{\nu}}{\infty} \leq \delta$, we have
			\begin{equation}\label{eq_lip_error}
				\sup_{\mb{\theta},\mb{\nu}} \vert J_r(u_{\mb{\nu}}) - J_r(u_{\mb{\theta}})\vert + \sup_{\mb{\theta},\mb{\nu}} \vert J_{r,N}(u_{\mb{\nu}}) - J_{r,N}(u_{\mb{\theta}})\vert \leq 2 \mathfrak{L} \norm{\mb{\theta} - \mb{\nu}}{\infty} \leq \frac{\varepsilon^2}{2}.
			\end{equation}
			For each $1 \leq i \leq N_{\delta}$, noting that $J_r(u_{\mb{\theta}_N^{*}}) = J_r(u_{\mb{\theta}_N^{*}}) - J_r(u_{\mb{\theta}_i}) + J_r(u_{\mb{\theta}_i}) - J_{r,N}(u_{\mb{\theta}_i}) + J_{r,N}(u_{\mb{\theta}_i}) - J_{r,N}(u_{\mb{\theta}_N^{*}}) + J_{r,N}(u_{\mb{\theta}_N^{*}})$, it follows that 
			\begin{equation}\label{eq_error_decompose}
				J_r(u_{\mb{\theta}_N^{*}}) \leq \vert J_r(u_{\mb{\theta}_N^{*}}) - J_r(u_{\mb{\theta}_i}) \vert + \vert J_r(u_{\mb{\theta}_i}) - J_{r,N}(u_{\mb{\theta}_i}) \vert + \vert J_{r,N}(u_{\mb{\theta}_i}) - J_{r,N}(u_{\mb{\theta}_N^{*}}) \vert + J_{r,N}(u_{\mb{\theta}_N^{*}}).
			\end{equation}
			Next, the infinite set $\Theta$ of trainable parameters is discretized by the $\delta$-net, then we use the estimate of residual for the $\delta$-net and combine it with the union bound to give the final estimate. Let $\mathcal{P}: \Theta \mapsto \bar{\mathsf{N}}$ be a projection (in the $\ell_{\infty}$ sense) onto $\bar{\mathsf{N}}$, i.e., $\mathcal{P}(\mb{\theta}) = \bar{\mb{\theta}}$ where $\bar{\mb{\theta}} = \argmin_{\mb{\upsilon} \in \bar{\mathsf{N}}} \norm{\mb{\theta} - \mb{\upsilon}}{\infty}$. Consider the following events for $1 \leq i \leq N_{\delta}$:
			\begin{equation*}
				\begin{aligned}
					\mathit{E_1} &= \{J_r(u_{\mb{\theta}_N^{*}}) \leq \varepsilon^2 + J_{r,N}(u_{\mb{\theta}_N^{*}})  \}, \\
					\mathit{E_{2,i}} &= \{J_r(u_{\mb{\theta}_i}) \leq \frac{\varepsilon^2}{2} + J_{r,N}(u_{\mb{\theta}_i})  \}, \\
					\mathit{E_{3,i}} &= \{\mathcal{P}(\mb{\theta}_N^{*}) = \mb{\theta}_i   \}, \\
					\mathit{E_4} &= \{ \exists \ i \in \{1, \ldots, N_{\delta}\}:  J_r(u_{\mb{\theta}_i}) \leq \frac{\varepsilon^2}{2} + J_{r,N}(u_{\mb{\theta}_i}) \ \text{and} \  \mathcal{P}(\mb{\theta}_N^{*}) = \mb{\theta}_i \}
				\end{aligned}
			\end{equation*}
			By \eqref{eq_lip_error} and \eqref{eq_error_decompose}, we known that if event $\mathit{E_4}$ occurs, then event $\mathit{E_1}$ occurs. Indeed, we have $\mathcal{P}(\mb{\theta}_N^{*}) = \mb{\theta}_i$ if $\mathit{E_4}$ occurs, which implies that $\norm{\mb{\theta}_N^{*} - \mb{\theta}_i}{\infty} \leq \delta$ due to the property of the $\delta$-net. By \eqref{eq_error_decompose} and using inequality \eqref{eq_lip_error} derived from the property of Lipschitz continuity, we can obtain $J_r(u_{\mb{\theta}_N^{*}}) \leq \varepsilon^2 + J_{r,N}(u_{\mb{\theta}_N^{*}})$, which means that $\mathit{E_1}$ occurs. Hence, we have $\mathit{E_4} \subseteq  \mathit{E_1}$, implying that 
			\begin{equation}\label{eq_prob_ineq}
				\mathbb{P}(\mathit{E_4}) \leq \mathbb{P}(\mathit{E_1}).
			\end{equation}
			
			According to the definition of $\delta$-net, we have 
			\begin{equation}\label{eq_sumeq1}
				\sum_i \mathbb{P}(\mathit{E_{3,i}}) = 1.
			\end{equation}
			By Assumption \ref{assump_bounded} and the Hoeffding inequality, we obtain
			\begin{equation}\label{eq_hoeff_ineq}
				\mathbb{P}(\mathit{E_{2,i}}) \geq 1 - \mathrm{exp}(\frac{-N_r \varepsilon^4}{2 c^2}).
			\end{equation}
			Noting that $\mathit{E_4} = \cup_{i=1}^{N_{\delta}} (\mathit{E_{2,i}} \cap \mathit{E_{3,i}})$, and combining \eqref{eq_prob_ineq}, \eqref{eq_sumeq1} and \eqref{eq_hoeff_ineq}, we have 
			\begin{equation*}
				\begin{aligned}
					\mathbb{P}(\mathit{E_1}) & \geq \mathbb{P}(\mathit{E_4}) = \sum\limits_i^{N_{\delta}} \mathbb{P}(\mathit{E_{2,i}} \cap \mathit{E_{3,i}}) = \sum\limits_{i=1}^{N_{\delta}} \left(\mathbb{P}(\mathit{E_{2,i}}) + \mathbb{P}(\mathit{E_{3,i}}) - \mathbb{P}( \mathit{E_{2,i}} \cup \mathit{E_{3,i}})  \right) \\
					&\geq 1 + \sum\limits_{i=1}^{N_{\delta}} \left( \mathbb{P}(\mathit{E_{2,i}}) - 1 \right) \\
					&\geq 1 - N_{\delta} \mathrm{exp}(\frac{-N_r \varepsilon^4}{2 c^2}) \\
					&\geq 1 - (4a \mathfrak{L}/\varepsilon^2)^D \mathrm{exp}(\frac{-N_r \varepsilon^4}{2 c^2}),
				\end{aligned}
			\end{equation*}
			which gives that
			\begin{equation*}
				J_r(u_{\mb{\theta}_N^{*}}) \leq \varepsilon^2 + J_{r,N}(u_{\mb{\theta}_N^{*}})
			\end{equation*}
			with probability at least $1 - (4a \mathfrak{L}/\varepsilon^2)^D \mathrm{exp}(-N_r \varepsilon^4/2c^2 )$.
		\end{proof}
		
		\section{Proof of Theorem \ref{thm_error_beha}}
		\begin{proof}
			Noting that 
			\begin{equation*}
				\mb{\theta}_N^{*, (k+1)} = \arg \min_{\mb{\theta}} \frac{1}{N_r} \sum_{i=1}^{N_r} \frac{r^2(\mb{x}^{(i)},\mb{\xi}^{(i)};\mb{\theta})}{p_{\mathsf{KRnet}}(\mb{x}^{(i)},\mb{\xi}^{(i)};\mb{\theta}_f^{*, (k)})}.
			\end{equation*} 
			Since $\mb{\theta}_N^{*, (k+1)}$ is the optimal solution at the $(k+1)$-th stage and $\mb{\theta}_{N}^{*,(k)}$ is used for initialization, we can obtain
			\begin{equation}\label{eq_discrete_rk}
				J_{r,N}(u_{\mb{\theta}_N^{*,(k+1)}}) = \frac{1}{N_r} \sum_{i=1}^{N_r} \frac{r^2(\mb{x}^{(i)},\mb{\xi}^{(i)};\mb{\theta}_N^{*,(k+1)})}{p_{\mathsf{KRnet}}(\mb{x}^{(i)},\mb{\xi}^{(i)};\mb{\theta}_f^{*, (k)})} \leq \frac{1}{N_r} \sum\limits_{i=1}^{N_r} \frac{r^2(\mb{x}^{(i)},\mb{\xi}^{(i)};\mb{\theta}_N^{*,(k)})}{p_{\mathsf{KRnet}}(\mb{x}^{(i)},\mb{\xi}^{(i)};\mb{\theta}_f^{*, (k)})}.
			\end{equation}
			Plugging $ p_{\mathsf{KRnet}}(\bx,\mb{\xi}; \btheta_f^{*,(k)}) = c_k r^2(\bx,\mb{\xi}; \btheta_N^{*,(k)})$ into \eqref{eq_discrete_rk}, we have
			\begin{equation*}
				J_{r,N}(u_{\mb{\theta}_N^{*,(k+1)}}) \leq \frac{1}{c_k}.
			\end{equation*}	
			Noting that $J_{r,N}(u_{\mb{\theta}_N^{*,(k+1)}})$ is a random variable and taking its expectation, it follows that
			\begin{equation*}
				\mathbb{E}(J_{r,N}(u_{\mb{\theta}_N^{*,(k+1)}}))  \leq \frac{1}{c_k} =  \int_{\Omega} r^2(\mb{x},\mb{\xi};\mb{\theta}_N^{*,(k)})d\mb{x}d\mb{\xi} = \mathbb{E}(J_{r,N}(u_{\mb{\theta}_N^{*,(k)}})),
			\end{equation*} 
			which completes the proof.
		\end{proof}
	\end{appendices}
	%
	%
   \bibliographystyle{elsarticle-harv} 
	\bibliography{tang}
	






\end{document}